\documentclass[12pt,draft]{article}
\usepackage{bm}
\usepackage{amsmath}
\usepackage{setspace}
\usepackage{amssymb}
\makeatletter
\newtheorem{lem}{Lemma}[section]
\newtheorem{prop}{Proposition}[section]
\newtheorem{thm}{Theorem}
\newtheorem{cor}{Corollary}[section]
\newtheorem{rem}{Remark}[section]
\newtheorem{exa}{Example}
\newtheorem{ass}{Assumption}

\usepackage{bm, txfonts}

\numberwithin{equation}{section}

\makeatother

\begin{document}

\title{The Bismut-Elworthy-Li type formulae for stochastic differential
equations with jumps}

\author{Atsushi TAKEUCHI%
\thanks{E-mail: \texttt{takeuchi@sci.osaka-cu.ac.jp} \protect \\
\quad \ \ Postal Address: Department of Mathematics, Osaka City
University, Sugimoto 3-3-138, Sumiyoshi-ku, Osaka 558-8585, JAPAN%
}}

\maketitle
{\allowdisplaybreaks{
\begin{abstract}
Consider jump-type stochastic differential equations with the drift,
diffusion and jump terms. Logarithmic derivatives of densities for
the solution process are studied, and the Bismut-Elworthy-Li type
formulae can be obtained under the uniformly elliptic condition on
the coefficients of the diffusion and jump terms. Our approach is
based upon the Kolmogorov backward equation by making full use of
the Markovian property of the process. 
\end{abstract}

\section{Introduction}

The Malliavin calculus has played an important role in many fields,
as one of powerful tools in infinite dimensional analysis. That has
also given us an attractive solution to the hypoelliptic problem for
the differential operator associated with a stochastic differential
equation, by means of probabilistic methods. It is well known that
the H\"ormander condition on the coefficients of the equation, which
is the condition about the Lie algebra generated by the vector fields
associated with the coefficients, yields the existence of the smooth
density function. See \cite{Bismut-1,Nualart} and references therein.
Bismut \cite{Bismut-2} also studied the logarithmic derivatives of
the density function with respect to the initial point of a stochastic
differential equation on Riemannian manifolds. His approach is based
upon the Girsanov transform on Brownian motions. The formula has a
nice flavour with the precise estimate of heat kernels or large deviation
principles. Elworthy and Li \cite{Elworthy and Li} also tackled the
same problem in more general class of stochastic differential equations
on Riemannian manifolds, via the martingale methods. Nowadays, the
celebrated formulae are called \textit{the Bismut-Elworthy-Li formulae}
after their great contributions. The logarithmic derivatives of the
density function is equivalent to the Greeks computations for pay-off
functions in mathematical finance. Fourni\'e \textit{et al.} \cite{Fournie et al-1}
applied the Malliavin calculus on the Wiener space to the sensitivity
analysis for asset price dynamics models. They also applied their
results to the numerical computations of the Greeks. 

All works stated above, paid attention to the case of the processes
without any jumps. There has been a natural and non-trivial question
whether a similar approach is applicable to the sensitivity analysis
in case of jump processes. The interests in jump processes are recently
getting more and more in mathematical finance. In the present paper,
we shall study the Bismut-Elworthy-Li type formulae for jump processes,
with respect to the initial point and the parameter governing the
equation. There are some approaches to tackle the problem on the sensitivities:
the Girsanov transform approach (\cite{Kawai-Takeuchi}) for L\'evy
processes initiated by Bismut \cite{Bismut-1}, the martingale methods
(\cite{Cass and Friz}) similarly to \cite{Elworthy and Li} in case
of diffusion processes, and an application of the Malliavin calculus
on the Wiener-Poisson space (\cite{Bally et al,Davis and Johansson,El-Khatib and Privault}).
In particular, Davis and Johansson \cite{Davis and Johansson}, and
Cass and Friz \cite{Cass and Friz} studied in case of jump diffusion
processes, but their approach does not take any effects from the jump
term. Bally \textit{et al.} \cite{Bally et al} studied the Malliavin
caluclus with respect to the jump amplitudes and the jump times, and
used the integration by parts formula in order to give numerical algorithms
for the sensitivity computations in a model driven by L\'evy processes.
The goal in the present paper is to compute the logarithmic derivatives
of densities, including not only the effect from the diffusion terms,
but also the one from the jump terms. The results obtained in this
paper corresponds to give another approach on the logarithmic derivatives
of the density studied by Bismut \cite{Bismut-1} using the Girsanov
transforms. 

This paper is organized as follows: In Section \ref{sec:Preliminaries},
we shall prepare some notations and introduce our framework. In Section
\ref{sec:Main-theorems}, the main theorems on the logarithmic derivatives
for the density with respect to the initial point and the parameter,
are given. Those proofs are done in Section \ref{sec:Proofs}. Some
typical examples are given in the final section.

\section{Preliminaries\label{sec:Preliminaries}}

At the beginning, we shall introduce some general notations. For $\alpha,\,\beta\in\mathbb{N}$,
denote by $C^{k}\left(\mathbb{R}^{\alpha}\,;\,\mathbb{R}^{\beta}\right)$
the class of $k$ times continuously differentiable, $\mathbb{R}^{\beta}$-valued
mappings on $\mathbb{R}^{\alpha}$, and by $C_{1+,b}^{k}\left(\mathbb{R}^{\alpha}\,;\,\mathbb{R}^{\beta}\right)$
the class of $C^{k}\left(\mathbb{R}^{\alpha}\,;\,\mathbb{R}^{\beta}\right)$-functions
with bounded derivatives of all orders more than 1. The subscript
$K$ of $C_{K}^{k}\left(\mathbb{R}^{\alpha}\,;\,\mathbb{R}^{\beta}\right)$
indicates the compact support. Denote by $\nabla=\left(\nabla_{1},\,\ldots,\,\nabla_{d}\right)$
the gradient operator in $\mathbb{R}^{d}$, by $\partial_{z}=\left(\partial_{z_{1}},\,\ldots,\,\partial_{z_{m}}\right)$
the one in $\mathbb{R}^{m}$, and by $\partial_{\varepsilon}=\left(\partial_{\varepsilon_{1}},\,\ldots,\,\partial_{\varepsilon_{l}}\right)$
the one in $\mathbb{R}^{l}$. Write the idensity by $I_{d}=\left(\delta_{jk}\,;\,1\le j,k\le d\right)\in\mathbb{R}^{d}\otimes\mathbb{R}^{d}$.
For $M\in\mathbb{R}^{\alpha}\otimes\mathbb{R}^{\beta}$, the symbol
$\left[M\right]_{ij}$ indicates the $\left(i,j\right)$-component
of $M$. For a subset $N\subset\mathbb{R}^{\alpha}$, denote its closure
by $\overline{N}$, its boundary by $\partial N$, and its Lebesgue
measure by $\left|N\right|$. For $\mathbb{R}^{d}\otimes\mathbb{R}^{m}$-valued
function $\Phi$ on $\mathbb{R}_{0}^{m}$ and $L\in\mathbb{R}^{m}\otimes\mathbb{R}^{l}\otimes\mathbb{R}^{m}$,
define\begin{gather*}
\left\{ \mbox{div}_{z}\left[\Phi\left(z\right)\right]\right\} _{k}=\sum_{i=1}^{m}\partial_{z_{i}}\left(\Phi_{ik}\ \left(z\right)\right),\ \mbox{div}_{z}\left[\Phi\left(z\right)\right]=\left(\left\{ \mbox{div}_{z}\left[\Phi\left(z\right)\right]\right\} _{1},\,\ldots,\,\left\{ \mbox{div}_{z}\left[\Phi\left(z\right)\right]\right\} _{d}\right),\\
\left\{ \mbox{Tr}\left[L\right]\right\} _{k}=\sum_{i=1}^{m}L_{iki},\quad\mbox{Tr}\left[L\right]=\left(\left\{ \mbox{Tr}\left[L\right]\right\} _{1},\,\ldots,\,\left\{ \mbox{Tr}\left[L\right]\right\} _{l}\right).\end{gather*}
Denote by $c_{i}$'s different positive finite constants. 

Write $\mathbb{R}_{0}^{m}=\mathbb{R}^{m}\backslash\left\{ \bm{0}\right\} $,
and let $d\nu$ be a L\'evy measure on $\mathbb{R}_{0}^{m}$. Moreover,
suppose that 

\begin{ass}\label{Levy measure} The measure $d\nu$ satisfies the
following three conditions: 
\begin{enumerate}
\item [(i)] for any $p\ge1$, \[
\int_{\mathbb{R}_{0}^{m}}\left(\left|z\right|I_{\left(\left|z\right|\le1\right)}+\left|z\right|^{p}I_{\left(\left|z\right|>1\right)}\right)d\nu<+\infty,\]

\item [(ii)] there exists a constant $\alpha>0$ such that, for any $\theta\in S^{m-1}$,
\[
\liminf_{\rho\searrow0}\rho^{\alpha}\int_{\mathbb{R}_{0}^{m}}\left(\left|z\cdot\theta/\rho\right|^{2}\wedge1\right)d\nu>0,\]

\item [(iii)] there exists a $C^{1}$-density $g\left(z\right)$ with respect
to the Lebesgue measure on $\mathbb{R}_{0}^{m}$ such that \[
\lim_{\left|z\right|\to\infty}|z|^{2}\, g(z)=0.\]

\end{enumerate}
\end{ass}

\begin{rem}{\rm{L\'evy processes such as tempered stable processes,
inverse Gaussian processes, etc. satsify Assumption \ref{Levy measure}.
\hfill{}$\square$}}\end{rem}

\begin{rem}{\rm{ In order to study the existence of a smooth density
for jump processes, Ishikawa and Kunita \cite{Ishikawa-Kunita}, and
Picard \cite{Picard} impose the following two conditions to the measure
$d\nu$ instead of Assumption \ref{Levy measure}: 
\begin{enumerate}
\item [{\it{(iv)}}] {\it{ there exists $0<\beta<2$ such that}} \[
\liminf_{\rho\searrow0}\rho^{-\beta}\int_{\left|z\right|\le\rho}\left|z\right|^{2}d\nu>0,\]

\item [{\it{(v)}}]  {\it{ there exists a positive definite matrix $\Xi\in\mathbb{R}^{m}\otimes\mathbb{R}^{m}$
such that, for any $\theta\in S^{m-1}$,}} \[
\liminf_{\rho\searrow0}\left(\int_{\left|z\right|\le\rho}\left|z\right|^{2}d\nu\right)^{-1}\int_{\left|z\right|\le\rho}\left|z\cdot\theta\right|^{2}d\nu\ge\theta\cdot\Xi\theta.\]

\end{enumerate}
The condition (iv) is called \textit{the order condition} on the measure
$d\nu$, and the L\'evy process satisfying the condition (v) is called
\textit{non-degenerate}. It can be easily checked that the above conditions
(iv) and (v) imply (ii) in Assumption \ref{Levy measure}. In fact,
\[
\int_{\mathbb{R}_{0}^{m}}\left(\left|z\cdot\theta/\rho\right|^{2}\wedge1\right)d\nu\ge\rho^{-2}\int_{\left|z\right|\le\rho}\left|z\cdot\theta\right|^{2}d\nu\ge c_{1}\,\rho^{-2+\beta}\,\theta\cdot\Xi\theta.\]
\hfill{}$\square$}}\end{rem}

Let $T>0$ be fixed, and $\left(\Omega,\mathcal{F},\mathbb{P}\right)$
the underlying probability space. Denote an $m$-dimensional Brownian
motion with $W_{0}=\bm{0}$ by $\left\{ W_{t}=\left(W_{t}^{1},\,\ldots,\, W_{t}^{m}\right)^{\ast}\,;\, t\in\left[0,T\right]\right\} $,
and by $d\mu$ the Poisson random measure on $\left[0,T\right]\times\mathbb{R}_{0}^{m}$
with the intensity $d\hat{\mu}=dt\, d\nu$. Let $\left\{ \mathcal{F}_{t}\,;\, t\in\left[0,T\right]\right\} $
be the augmented filtration generated by $W$ and $d\mu$ with respect
to $\mathbb{P}$. Define $d\tilde{\mu}=d\mu-d\hat{\mu}$ and $d\overline{\mu}=I_{\left(\left|z\right|\le1\right)}d\tilde{\mu}+I_{\left(\left|z\right|>1\right)}d\mu$.
Let $a_{i}:\mathbb{R}^{l}\times\mathbb{R}^{d}\to\mathbb{R}^{d}\ \left(i=0,\,1,\,\ldots,\, m\right)$
and $b:\mathbb{R}^{l}\times\mathbb{R}^{d}\times\mathbb{R}_{0}^{m}\to\mathbb{R}^{d}$
such that 

\begin{ass}\label{assumption: coefficients} The $\mathbb{R}^{d}$-valued
functions $a_{0},\, a_{1},\,\ldots,\, a_{m},\, b$ satisfy 
\begin{enumerate}
\item [(i)] $a_{i}\left(\varepsilon,\ \cdot\ \right)\in C_{1+,b}^{\infty}\left(\mathbb{R}^{d}\,;\,\mathbb{R}^{d}\right)$
for each $\varepsilon\in\mathbb{R}^{l}$, and $a_{i}\left(\ \cdot\ ,y\right)\in C_{1+,b}^{1}\left(\mathbb{R}^{l}\,;\,\mathbb{R}^{d}\right)$
for each $y\in\mathbb{R}^{d}$, 
\item [(ii)] $b_{\cdot}\left(\varepsilon,\ \cdot\ \right)\in C_{1+,b}^{\infty,\infty}\left(\mathbb{R}^{d}\times\mathbb{R}_{0}^{m}\,;\,\mathbb{R}^{d}\right)$
for each $\varepsilon\in\mathbb{R}^{l}$, and $b_{z}\left(\ \cdot\ ,y\right)\in C_{1+,b}^{1}\left(\mathbb{R}^{l}\,;\,\mathbb{R}^{d}\right)$
for each $\left(y,z\right)\in\mathbb{R}^{d}\times\mathbb{R}_{0}^{m}$, 
\item [(iii)] for each $\varepsilon\in\mathbb{R}^{l}$, \[
\lim_{\left|z\right|\searrow0}b_{z}\left(\varepsilon,y\right)=0,\quad\inf_{y\in\mathbb{R}^{d}}\inf_{z\in\mathbb{R}_{0}^{m}}\ \left|\det\left[I_{d}+\nabla b_{z}\left(\varepsilon,y\right)\right]\right|>0.\]

\end{enumerate}
\end{ass}

\medskip{}

Write $a=\left(a_{1},\,\ldots,\, a_{m}\right)$. For $\left(\varepsilon,x\right)\in\mathbb{R}^{l}\times\mathbb{R}^{d}$,
consider the stochastic differential equation: \begin{equation}
{\displaystyle dx_{t}=a_{0}\left(\varepsilon,x_{t}\right)dt+a\left(\varepsilon,x_{t}\right)\circ dW_{t}{\displaystyle +\int_{\mathbb{R}_{0}^{m}}b_{z}\left(\varepsilon,x_{t-}\right)d\overline{\mu}},}\quad x_{0}=x.\label{eq:SDE}\end{equation}
Since the coefficients satisfy the Lipschitz and linear growth conditions
under Assumption \ref{assumption: coefficients}, there exists a unique
solution $\left\{ x_{t}=x_{t}^{x,\varepsilon}\,;\, t\in\left[0,T\right]\right\} $
(cf. \cite{Ikeda-Watanabe}). The infinitesimal generator $\mathcal{L}^{\varepsilon}$
associated with the solution process $\left\{ x_{t}\,;\, t\in\left[0,T\right]\right\} $
is given by\[
\begin{split}\left(\mathcal{L}^{\varepsilon}f\right)\left(y\right) & =A_{0}^{\varepsilon}f\left(y\right)+\frac{1}{2}\sum_{i=1}^{m}A_{i}^{\varepsilon}A_{i}^{\varepsilon}f\left(y\right)\ +\int_{\mathbb{R}_{0}^{m}}\left\{ \mathfrak{B}_{z}^{\varepsilon}f\left(y\right)-B_{z}^{\varepsilon}f\left(y\right)I_{\left(|z|\le1\right)}\right\} d\nu\end{split}
\]
for $f\in C_{K}^{2}\left(\mathbb{R}^{d};\mathbb{R}\right)$, where
$A_{i}^{\varepsilon}f\left(y\right)=\nabla f\left(y\right)a_{i}\left(\varepsilon,y\right)$,
$B_{z}^{\varepsilon}f\left(y\right)=\nabla f\left(y\right)b_{z}\left(\varepsilon,y\right)$
are vector fields, $A_{i}^{\varepsilon}A_{i}^{\varepsilon}f\left(y\right)=\nabla\left(\nabla f\left(y\right)\, a_{i}\left(\varepsilon,y\right)\right)\, a_{i}\left(\varepsilon,y\right)$
and $\mathfrak{B}_{z}^{\varepsilon}f\left(y\right)=f\left(y+b_{z}\left(\varepsilon,y\right)\right)-f\left(y\right)$.
Moreover, Assumption \ref{assumption: coefficients} yields that the
mapping $\mathbb{R}^{d}\ni x\mapsto x_{t}\in\mathbb{R}^{d}$ has a
$C^{1}$-modification for each $\left(t,\varepsilon\right)\in\left[0,T\right]\times\mathbb{R}^{l}$,
and its Jacobi matrix $Z_{t}:=\nabla_{x}x_{t}$ satisfies the linear
stochastic differential equation: \begin{equation}
{\displaystyle {\displaystyle dZ_{t}=\nabla a_{0}\left(\varepsilon,x_{t}\right)Z_{t}dt+\nabla a\left(\varepsilon,x_{t}\right)Z_{t}\circ dW_{t}}{\displaystyle +\int_{\mathbb{R}_{0}^{m}}\nabla b_{z}\left(\varepsilon,x_{t-}\right)Z_{t-}d\overline{\mu},\quad Z_{0}=I_{d}.}}\label{eq:Jacobi matrix}\end{equation}
Let $\left\{ U_{t}\,;\, t\in\left[0,T\right]\right\} $ be the solution
to the linear stochastic differential equation: $U_{0}=I_{d}$ and
\begin{equation}
\begin{split}dU_{t} & =-U_{t}\,\nabla a_{0}\left(\varepsilon,x_{t}\right)dt-U_{t}\,\nabla a\left(\varepsilon,x_{t}\right)\circ dW_{t}{\displaystyle {\displaystyle -\int_{\mathbb{R}_{0}^{m}}U_{t-}\left[\left(I_{d}+\nabla b_{z}\right)^{-1}\nabla b_{z}\right]\left(\varepsilon,x_{t-}\right)d\overline{\mu}}}\\
 & \quad+{\displaystyle \int_{|z|\le1}U_{t}\,\left[\left(I_{d}+\nabla b_{z}\right)^{-1}\left(\nabla b_{z}\right)^{2}\right]\left(\varepsilon,x_{t}\right)d\hat{\mu}.}\end{split}
\label{eq:inverse matrix}\end{equation}
Then, $Z_{t}\, U_{t}=U_{t}\, Z_{t}=I_{d}$ holds for each $t\in\left[0,T\right]$
by the It\^o product formula. Under Assumption \ref{Levy measure}
on the measure $d\nu$, and Assumption \ref{assumption: coefficients}
on the coefficients, the upper estimate \[
\mathbb{E}\left[\sup_{t\in\left[0,T\right]}\left\{ \left|x_{t}\right|^{p}+\left\Vert Z_{t}\right\Vert ^{p}+\left\Vert U_{t}\right\Vert ^{p}\right\} \right]\le c_{2,p,\varepsilon,T}\left(1+\left|x\right|^{p}\right)\]
holds for any $p>1$. See \cite{Fujiwara-Kunita}. Moreover, we have 

\begin{prop}\label{differentiability in varepsilon} For each $t\in\left[0,T\right]$,
the mapping $\mathbb{R}^{l}\ni\varepsilon\mapsto x_{t}\in\mathbb{R}^{d}$
has a $C^{1}$-modification, and the derivative $H_{t}:=\partial_{\varepsilon}x_{t}$
satisfies the equation: $H_{0}=\bm{0}\in\mathbb{R}^{l}\otimes\mathbb{R}^{d}$,
and \begin{equation}
\begin{aligned}dH_{t} & =\nabla a_{0}\left(\varepsilon,x_{t}\right)\, H_{t}\, dt+\nabla a\left(\varepsilon,x_{t}\right)\, H_{t}\circ dW_{t}+\int_{\mathbb{R}_{0}^{m}}\nabla b_{z}\left(\varepsilon,x_{t-}\right)\, H_{t-}\, d\overline{\mu}\\
 & \quad+\partial_{\varepsilon}a_{0}\left(\varepsilon,x_{t}\right)\, dt+\partial_{\varepsilon}a\left(\varepsilon,x_{t}\right)\circ dW_{t}+\int_{\mathbb{R}_{0}^{m}}\partial_{\varepsilon}b_{z}\left(\varepsilon,x_{t-}\right)\, d\overline{\mu}.\end{aligned}
\label{eq:partial_epsilon x_t}\end{equation}
Moreover, it holds that, for any $p>1$ and any compact subset $K$
in $\mathbb{R}^{l}$, \[
\sup_{\varepsilon\in K}\mathbb{E}\left[\sup_{t\in\left[0,T\right]}\left\Vert H_{t}\right\Vert ^{p}\right]\le c_{3,p,T}.\]
\end{prop}

\textit{Proof}. We shall write $x_{t}=x_{t}^{\varepsilon}$, in order
to emphasize the dependence on $\varepsilon\in\mathbb{R}^{l}$ throughout
the proof. Let $\left(\varepsilon,\delta\right)\in\mathbb{R}^{l}\times\mathbb{R}^{l}$.
Since \begin{align*}
x_{t}^{\varepsilon}-x_{t}^{\delta} & =\int_{0}^{t}\left\{ a_{0}\left(\varepsilon,x_{s}^{\varepsilon}\right)-a_{0}\left(\delta,x_{s}^{\delta}\right)\right\} \, ds+\int_{0}^{t}\left\{ a\left(\varepsilon,x_{s}^{\varepsilon}\right)-a\left(\delta,x_{s}^{\delta}\right)\right\} \circ dW_{s}\\
 & \qquad+\int_{0}^{t}\!\!\int_{\mathbb{R}_{0}^{m}}\left\{ b_{z}\left(\varepsilon,x_{s-}^{\varepsilon}\right)-b_{z}\left(\delta,x_{s-}^{\delta}\right)\right\} \, d\overline{\mu},\end{align*}
we can get \[
\mathbb{E}\left[\sup_{t\le T}\left|x_{t}^{\varepsilon}-x_{t}^{\delta}\right|^{p}\right]\le c_{4,p,x,T}\,\left|\varepsilon-\delta\right|^{p}\]
for any $p>1$, from Assumption \ref{Levy measure} and \ref{assumption: coefficients}.
Thus, the Kolmogorov continuity criterion tells us that the mapping
$\mathbb{R}^{l}\ni\varepsilon\longmapsto x_{t}^{\varepsilon}\in\mathbb{R}^{d}$
has a continuous modification for each $t\ge0$ and $x\in\mathbb{R}^{d}$. 

Next, we shall study the differentiability of $x_{t}^{\varepsilon}$
in $\varepsilon\in\mathbb{R}^{l}$. Let $0\neq\xi,\,\zeta\in\mathbb{R}$,
and $\bm{e}_{k}\in\mathbb{R}^{l}$ the $k$-th unit vector. Since
\begin{align*}
\left(x_{t}^{\varepsilon+\xi\bm{e}_{k}}-x_{t}^{\varepsilon}\right)/\xi & =\int_{0}^{t}\left\{ a_{0}\left(\varepsilon+\xi\bm{e}_{k},x_{s}^{\varepsilon+\xi\bm{e}_{k}}\right)-a_{0}\left(\varepsilon,x_{s}^{\varepsilon}\right)\right\} /\xi\, ds\\
 & \qquad+\int_{0}^{t}\left\{ a\left(\varepsilon+\xi\bm{e}_{k},x_{s}^{\varepsilon+\xi\bm{e}_{k}}\right)-a\left(\varepsilon,x_{s}^{\varepsilon}\right)\right\} /\xi\circ dW_{s}\\
 & \qquad+\int_{0}^{t}\!\!\int_{\mathbb{R}_{0}^{m}}\left\{ b_{z}\left(\varepsilon+\xi\bm{e}_{k},x_{s-}^{\varepsilon+\xi\bm{e}_{k}}\right)-b_{z}\left(\varepsilon,x_{s-}^{\varepsilon}\right)\right\} /\xi\, d\overline{\mu},\end{align*}
we can get the upper estimate\[
\mathbb{E}\left[\sup_{t\le T}\left|\left(x_{t}^{\varepsilon+\xi\bm{e}_{k}}-x_{t}^{\varepsilon}\right)/\xi-\left(x_{t}^{\varepsilon+\zeta\bm{e}_{k}}-x_{t}^{\varepsilon}\right)/\zeta\right|^{p}\right]\le c_{5,p,x,T,\varepsilon,k}\,|\xi-\zeta|^{p}\]
for any $p>1$. Hence, the mapping $\mathbb{R}^{l}\ni\varepsilon\longmapsto x_{t}^{\varepsilon}\in\mathbb{R}^{d}$
has a $C^{1}$-modification with respect to the parameter $\varepsilon\in\mathbb{R}^{l}$
for each $t\ge0$ and $x\in\mathbb{R}^{d}$, via the Kolmogorov continuity
criterion, again. 

Furthermore, Assumption \ref{assumption: coefficients} enables us
to justify that the derivative $\partial_{\varepsilon}x_{t}^{\varepsilon}$
satisfies the equation \eqref{eq:partial_epsilon x_t}. It is an easy
work to check the upper estimate of $\partial_{\varepsilon}x_{t}$
in the assertion. \hfill{}$\square$ 

\smallskip{}

\begin{cor}\label{Corollary: vega} The derivative $H_{t}=\partial_{\varepsilon}x_{t}$
can be computed as follows: \begin{equation}
\begin{aligned}H_{t} & =Z_{t}\int_{0}^{t}U_{s}\left(\partial_{\varepsilon}a_{0}\left(\varepsilon,x_{s}\right)-\int_{\left|z\right|\le1}\left[\left(I_{d}+\nabla b_{z}\right)^{-1}\nabla b_{z}\,\partial_{\varepsilon}b_{z}\right]\left(\varepsilon,x_{s}\right)d\nu\right)ds\\
 & \quad+Z_{t}\int_{0}^{t}U_{s}\partial_{\varepsilon}a\left(\varepsilon,x_{s}\right)\circ dW_{s}+Z_{t}\int_{0}^{t}\!\!\int_{\mathbb{R}_{0}^{m}}U_{s-}\left[\left(I_{d}+\nabla b_{z}\right)^{-1}\partial_{\varepsilon}b_{z}\right]\left(\varepsilon,x_{s-}\right)d\overline{\mu}\\
 & \left(=:Z_{t}\int_{0}^{t}f_{0}^{\varepsilon}\left(s\right)ds+Z_{t}\int_{0}^{t}f^{\varepsilon}\left(s\right)\circ dW_{s}+Z_{t}\int_{0}^{t}\!\!\int_{\mathbb{R}_{0}^{m}}h_{z}^{\varepsilon}\left(s\right)d\overline{\mu}.\right)\end{aligned}
\end{equation}
\end{cor}

\textit{Proof}. Obvious by Proposition \ref{differentiability in varepsilon}
and the It\^o product formula. \hfill{}$\square$

\section{Main theorems\label{sec:Main-theorems}}

Let us present an assumption on the coefficients of the equation \eqref{eq:SDE},
which is crucial for discussions in what follows. 

\begin{ass}\label{Uniform elliptic}There exist constants $c_{6,\varepsilon},\, c_{7,\varepsilon}>0$
such that \[
\sum_{i=1}^{m}\left|\xi\cdot a_{i}\left(\varepsilon,y\right)\right|^{2}\ge c_{6,\varepsilon}\left|\xi\right|^{2},\quad\sum_{i=1}^{m}\left|\xi\cdot\partial_{z_{i}}b_{z}\left(\varepsilon,y\right)\right|^{2}\ge c_{7,\varepsilon}\left|\xi\right|^{2}\]
for any $\left(y,\xi\right)\in\mathbb{R}^{d}\times\mathbb{R}^{d}$
and $z\in\mathbb{R}_{0}^{m}$. \end{ass}\bigskip{}

Define the $\mathbb{R}^{d}$-valued function $\tilde{b}$ by $\tilde{b}_{z}\left(\varepsilon,y\right)=\left[\left(I_{d}+\nabla b_{z}\right)^{-1}\partial_{z}b_{z}\right]\left(\varepsilon,y\right)\, z$.
Then, we shall introduce the well-known criterion on the existence
of the smooth density.  

\begin{prop}[cf. \cite{Komatsu-Takeuchi-1, Komatsu-Takeuchi-2}]\label{KT}
If there exist constants $c_{8,\varepsilon}>0$ and $\gamma>0$ such
that \begin{equation}
\inf_{y\in\mathbb{R}^{d}}\ \inf_{\xi\in S^{d-1}}\left\{ \sum_{i=1}^{m}\left|a_{i}\left(\varepsilon,y\right)\cdot\xi/\rho\right|^{2}+\int_{\mathbb{R}_{0}^{m}}\left(\left|\tilde{b}_{z}\left(\varepsilon,y\right)\cdot\xi/\rho\right|^{2}\wedge1\right)\, d\nu\right\} \ge c_{8,\varepsilon}\,\rho^{-\gamma}\label{uniformly elliptic condition}\end{equation}
for $0<\rho<1$, then the probability law of the random variable $x_{T}=x_{T}^{x,\varepsilon}$
has a density $p_{T}\left(\varepsilon,x,y\right)$ with respect to
the Lebesgue measure on $\mathbb{R}^{d}$ such that $p_{T}\left(\varepsilon,x,y\right)$
is smooth in $y\in\mathbb{R}^{d}$. \end{prop}

\begin{rem}\label{known fact}{\rm{ It can be easily checked that
Assumption \ref{Levy measure}, \ref{assumption: coefficients} and
\ref{Uniform elliptic} imply the condition \eqref{uniformly elliptic condition}
in Proposition \ref{KT}. In fact, since \begin{align*}
1 & \le\left|\left[\left(I_{d}+\nabla b_{z}\right)^{-1}\right]^{\ast}\left(\varepsilon,y\right)\,\xi\right|\left|\left[I_{d}+\nabla b_{z}\right]\left(\varepsilon,y\right)\,\xi\right|\le c_{9,\varepsilon}\,\left|\left[\left(I_{d}+\nabla b_{z}\right)^{-1}\right]^{\ast}\left(\varepsilon,y\right)\,\xi\right|\end{align*}
for $\xi\in S^{d-1}$, we see that \begin{align*}
\left|\left[\left(I_{d}+\nabla b_{z}\right)^{-1}\partial_{z}b_{z}\right]^{\ast}\left(\varepsilon,y\right)\,\xi/\rho\right|^{2} & \ge c_{10,\varepsilon}\,\left|\left[\left(I_{d}+\nabla b_{z}\right)^{-1}\right]^{\ast}\left(\varepsilon,y\right)\,\xi/\rho\right|^{2}\ge c_{11,\varepsilon}\,\rho^{-2}\end{align*}
 under Assumption \ref{assumption: coefficients} and \ref{Uniform elliptic}.
Then we have\begin{align*}
\int_{\mathbb{R}_{0}^{m}}\left(\left|\tilde{b}_{z}\left(\varepsilon,y\right)\cdot\xi/\rho\right|^{2}\wedge1\right)\, d\nu & =\int_{\mathbb{R}_{0}^{m}}\left(\left|\left(\left[\left(I_{d}+\nabla b_{z}\right)^{-1}\partial_{z}b_{z}\right]^{\ast}\left(\varepsilon,y\right)\,\xi/\rho\right)\cdot z\right|^{2}\wedge1\right)\, d\nu\\
 & \ge c_{12,\varepsilon}\,\inf_{\left|\theta\right|=1}\int_{\mathbb{R}_{0}^{m}}\left(\left|z\cdot\theta/\rho\right|^{2}\wedge1\right)\, d\nu\\
 & \ge c_{13,\varepsilon}\,\rho^{-\alpha}\end{align*}
from Assumption \ref{Levy measure} (ii), for sufficiently small $0<\rho<1$.
\hfill{}$\square$

}}\end{rem}\bigskip{}

We are now in a position to present main results. To avoid lengthy
expressions, let us prepare some auxiliary notations. Define \begin{gather*}
v_{s}\left(\varepsilon,z\right)=\left[\left(\partial_{z}b_{z}\right)^{-1}\left(I_{d}+\nabla b_{z}\right)\right]\left(\varepsilon,x_{s}\right)\, Z_{s}\,\left|z\right|^{2},\ \\
J_{\tau,t}=\int_{\tau}^{t}\!\!\int_{\mathbb{R}_{0}^{m}}\frac{\mbox{div}_{z}\left[g\left(z\right)v_{s}\left(\varepsilon,z\right)\right]}{g\left(z\right)}d\tilde{\mu},\quad K_{\tau,t}=\int_{\tau}^{t}\!\!\int_{\mathbb{R}_{0}^{m}}2z^{\ast}v_{s}\left(\varepsilon,z\right)d\mu,\ \\
A_{\tau,t}=\left(t-\tau\right)+\int_{\tau}^{t}\!\!\int_{\mathbb{R}_{0}^{m}}\left|z\right|^{2}d\mu,\quad L_{\tau,t}=\int_{\tau}^{t}dW_{s}^{\ast}\, a\left(\varepsilon,x_{s}\right)^{-1}Z_{s},\ \ \end{gather*}
for $0\le\tau\le t\le T$. We first derive the sensitivity formula
with respect to $x\in\mathbb{R}^{d}$. 

\begin{thm}\label{Theorem 1}Let $\varphi$ be in $C_{K}^{2}\left(\mathbb{R}^{d}\,;\,\mathbb{R}\right)$.
Then, it holds that \begin{equation}
\nabla_{x}\left(\mathbb{E}\left[\varphi\left(x_{T}\right)\right]\right)=\mathbb{E}\left[\varphi\left(x_{T}\right)\left\{ \frac{L_{0,T}-J_{0,T}}{A_{0,T}}+\frac{K_{0,T}}{\left(A_{0,T}\right)^{2}}\right\} \right]\left(=:\mathbb{E}\left[\varphi\left(x_{T}\right)\Gamma_{T}^{\left(1\right)}\right]\right).\label{eq:BEL 1st}\end{equation}
 \end{thm}\bigskip{}

Next, we shall study the sensitivity in $\varepsilon\in\mathbb{R}^{l}$.
For the sake of simplicity on notations, define \begin{gather*}
L_{t}^{\varepsilon}=\int_{0}^{t}dW_{s}^{\ast}\, a\left(\varepsilon,x_{s}\right)^{-1}Z_{s}\, f_{0}^{\varepsilon}\left(s\right),\quad G_{t}^{\varepsilon}=\int_{0}^{t}f^{\varepsilon}\left(s\right)\circ dW_{s},\\
R_{t}^{\varepsilon}=\frac{1}{t}\int_{0}^{t}dW_{s}^{\ast}\, a\left(\varepsilon,x_{s}\right)^{-1}Z_{s}\, G_{t}^{\varepsilon},\quad Q_{t}^{\varepsilon}=\frac{1}{t}\int_{0}^{t}\mbox{Tr}\left[a\left(\varepsilon,x_{s}\right)^{-1}Z_{s}\, D_{s}G_{t}^{\varepsilon}\right]ds,\ \\
\tilde{v}_{s}\left(\varepsilon,z\right)=\left[\left(\partial_{z}b_{z}\right)^{-1}\left(I_{d}+\nabla b_{z}\right)\right]\left(\varepsilon,x_{s}\right)\, Z_{s}\, h_{z}^{\varepsilon}\left(s\right),\quad J_{t}^{\varepsilon}=\int_{0}^{t}\!\!\int_{\mathbb{R}_{0}^{m}}\frac{\mbox{div}_{z}\left[g\left(z\right)\tilde{v}_{s}\left(\varepsilon,z\right)\right]}{g\left(z\right)}d\tilde{\mu},\end{gather*}
where $f_{0}^{\varepsilon}\left(s\right)$, $f^{\varepsilon}\left(s\right)$
and $h_{z}^{\varepsilon}\left(s\right)$ are given in Corollary \ref{Corollary: vega},
and $\left\{ D_{s}\,;\, s\in\left[0,T\right]\right\} $ is the Malliavin
derivative operator. 

\begin{thm}\label{vega}Let $\varphi\in C_{K}^{2}\left(\mathbb{R}^{d}\,;\,\mathbb{R}\right)$.
Then, it holds that \begin{equation}
\partial_{\varepsilon}\left(\mathbb{E}\left[\varphi\left(x_{T}\right)\right]\right)=\mathbb{E}\left[\varphi\left(x_{T}\right)\left\{ L_{T}^{\varepsilon}+R_{T}^{\varepsilon}-Q_{T}^{\varepsilon}-J_{T}^{\varepsilon}\right\} \right]\left(=:\mathbb{E}\left[\varphi\left(x_{T}\right)\Gamma_{T}^{\left(2\right)}\right]\right).\end{equation}
 \end{thm}

\bigskip{}

Finally, we next study the second order derivative in $x\in\mathbb{R}^{d}$.
To keep the presentation as concise as possible, write $\tilde{T}=T/2,$
and define \begin{gather*}
F_{ijk}^{1}\left(\varepsilon,t\right)=\sum_{\beta=1}^{d}\left\{ \nabla_{x_{j}}\left[a\left(\varepsilon,x_{t}\right)^{-1}\right]_{i\beta}Z_{t}^{\beta k}+\left[a\left(\varepsilon,x_{t}\right)^{-1}\right]_{i\beta}\nabla_{x_{k}}Z_{t}^{\beta j}\right\} ,\\
F_{ijk}^{2}\left(\varepsilon,t,z\right)=-\sum_{\beta=1}^{d}\left[\left\{ \partial_{z}b_{z}\left(\varepsilon,x_{t}\right)\right\} ^{-1}\right]_{i\beta}\nabla_{x_{k}}\left[\left(I_{d}+\nabla b_{z}\left(\varepsilon,x_{t}\right)\right)Z_{t}\right]_{\beta j}\,\left|z\right|^{2},\\
F_{ijk}^{3}\left(\varepsilon,t,z\right)=\sum_{\beta=1}^{d}\left[\left\{ \partial_{z}b_{z}\left(\varepsilon,x_{t}\right)\right\} ^{-1}\right]_{i\beta}\left[\partial_{z}\nabla_{x_{j}}\left(b_{z}\left(\varepsilon,x_{t}\right)\right)\, v_{t}\left(\varepsilon,z\right)\right]_{\beta jk}\end{gather*}
for $1\le i\le m$ and $1\le j,k\le d$. Moreover, write $F^{1}\left(\varepsilon,t\right)=\left(F_{ijk}^{1}\left(\varepsilon,t\right)\right)_{1\le i\le m,1\le j,k\le d}$,
and $F^{\sigma}\left(\varepsilon,t,z\right)=\left(F_{ijk}^{\sigma}\left(\varepsilon,t,z\right)\right)_{1\le i\le m,1\le j,k\le d}\ $
for $\sigma=2,\,3$. Then we have 

\begin{thm}\label{Theorem 2}Let $\varphi$ be in $C_{K}^{2}\left(\mathbb{R}^{d}\,;\,\mathbb{R}\right)$.
Then, the equality\begin{equation}
\begin{aligned} & \nabla_{x}\nabla_{x}\left(\mathbb{E}\left[\varphi\left(x_{T}\right)\right]\right)\\
 & =\mathbb{E}\Bigg[\varphi\left(x_{T}\right)\Bigg\{\left(\Gamma_{\tilde{T},T}^{\left(1\right)\ast}+\frac{K_{\tilde{T},T}^{\ast}}{\ A_{\tilde{T},T}\, A_{0,\tilde{T}}\ }\right)\Gamma_{0,\tilde{T}}^{\left(1\right)}+\frac{K_{\tilde{T},T}^{\ast}\, K_{0,\tilde{T}}}{A_{\tilde{T},T}\,\left(A_{0,\tilde{T}}\right)^{3}}-\sum_{\sigma=2}^{3}\int_{0}^{\tilde{T}}\!\!\int_{\mathbb{R}_{0}^{m}}\frac{\mbox{div}_{z}\left[F^{\sigma}\left(\varepsilon,s,z\right)\right]\ }{A_{0,\tilde{T}}+\left|z\right|^{2}}d\hat{\mu}\\
 & \qquad\qquad\quad+\frac{1}{A_{0,\tilde{T}}}\left(\int_{0}^{\tilde{T}}F^{1}\left(\varepsilon,s\right)dW_{s}+\sum_{\sigma=2}^{3}\int_{0}^{\tilde{T}}\!\!\int_{\mathbb{R}_{0}^{m}}\mbox{div}_{z}\left[F^{\sigma}\left(\varepsilon,s,z\right)\right]d\mu\right)\Bigg\}\Bigg]\\
 & \left(=:\mathbb{E}\left[\varphi\left(x_{T}\right)\Gamma_{T}^{\left(3\right)}\right]\right)\end{aligned}
\label{eq:BEL 2nd}\end{equation}
 holds, where $\Gamma_{\tau,t}^{\left(1\right)}=\left(L_{\tau,t}-J_{\tau,t}\right)/A_{\tau,t}+K_{\tau,t}/\left(A_{\tau,t}\right)^{2}$
for $0\le\tau\le t\le T$. 

\end{thm}

\begin{rem}\label{inverse integrability}\rm{ For $0\le\tau<t\le T$,
define $N_{\tau,t}^{\lambda}=\int_{\tau}^{t}\!\!\int_{\mathbb{R}_{0}^{m}}\left(e^{-\lambda\left|z\right|^{2}}-1\right)\, d\hat{\mu}$.
Then, $\mathbb{E}\left[A_{\tau,t}^{-p}\right]<+\infty$ holds for
any $p>1$, since the condition (ii) in Assumption \ref{Levy measure}
on the measure $d\nu$ yields that \begin{align*}
\mathbb{E}\left[A_{\tau,t}^{-p}\right] & =\frac{1}{\Gamma\left(p\right)}\int_{0}^{+\infty}\lambda^{p-1}\,\mathbb{E}\left[\exp\left(-\lambda A_{\tau,t}\right)\right]\, d\lambda\\
 & =\frac{1}{\Gamma\left(p\right)}\int_{0}^{+\infty}\lambda^{p-1}\,\mathbb{E}\left[\exp\left(-\lambda A_{\tau,t}-N_{\tau,t}^{\lambda}\right)\right]\, e^{N_{\tau,t}^{\lambda}}d\lambda\\
 & \le c_{14,T}\,\int_{0}^{+\infty}\lambda^{p-1}\,\exp\left\{ -\left(t-\tau\right)\lambda-\left(t-\tau\right)\int_{\mathbb{R}_{0}^{m}}\left\{ \left(\lambda\left|z\right|^{2}\right)\wedge1\right\} \, d\nu\right\} \, d\lambda\\
 & \le c_{14,T}\,\int_{0}^{+\infty}\lambda^{p-1}\,\exp\left\{ -\left(t-\tau\right)\lambda-c_{15}\left(t-\tau\right)\lambda^{\alpha/2}\right\} \, d\lambda\\
 & <+\infty.\end{align*}
\hfill{}$\square$}\end{rem}\medskip{}

Denote by $\mathcal{U}$ the family of bounded domains and their complements
in $\mathbb{R}^{d}$. Define the class $\mathfrak{F}$ of $\mathbb{R}$-valued
functions by \[
\mathfrak{F}=\left\{ f=\sum_{k=1}^{n}\alpha_{k}\, f_{k}\, I_{A_{k}}\,;\, n\in\mathbb{N},\,\alpha_{k}\in\mathbb{R},\,\left|f_{k}\left(y\right)\right|\le c_{16,k}\left(1+\left|y\right|\right),\, A_{k}\in\mathcal{U}\ \right\} .\]

\begin{cor}\label{sensitivity formulae} Let $\varphi\in\mathfrak{F}$,
and $\Gamma_{T}^{\left(i\right)}\ \left(i=1,\,2,\,3\right)$ be random
variables defined in Theorem \ref{Theorem 1}, \ref{vega} and \ref{Theorem 2}.
Then, the following equalities hold.\begin{gather*}
\nabla_{x}\left(\mathbb{E}\left[\varphi\left(x_{T}\right)\right]\right)=\mathbb{E}\left[\varphi\left(x_{T}\right)\Gamma_{T}^{\left(1\right)}\right],\ \\
\partial_{\varepsilon}\left(\mathbb{E}\left[\varphi\left(x_{T}\right)\right]\right)=\mathbb{E}\left[\varphi\left(x_{T}\right)\Gamma_{T}^{\left(2\right)}\right],\ \\
\nabla_{x}\nabla_{x}\left(\mathbb{E}\left[\varphi\left(x_{T}\right)\right]\right)=\mathbb{E}\left[\varphi\left(x_{T}\right)\Gamma_{T}^{\left(3\right)}\right].\end{gather*}
 \end{cor}

\begin{rem}{\rm{ The class $\mathfrak{F}$ is smaller than that
of measurable functions $\varphi$ satisfying $\mathbb{E}\left[\left|\varphi\left(x_{T}\right)\right|^{2}\right]<+\infty$.
But, the class $\mathfrak{F}$ is rich enough from a practical point
in mathematical finance, because various payoff functions for asset
price dynamics such as call options, put options, digital options,
and so on, are included in $\mathfrak{F}$. \hfill{}$\square$}}\end{rem}

\begin{rem}\label{uniformly elliptic diffusion}{\rm{ Consider the
case where the $\mathbb{R}^{d}\otimes\mathbb{R}^{d}$-valued function
$a\left(\varepsilon,y\right)a\left(\varepsilon,y\right)^{\ast}$ is
uniformly elliptic in $y\in\mathbb{R}^{d}$, while the function $\partial_{z}b_{z}\left(\varepsilon,y\right)\partial_{z}b_{z}\left(\varepsilon,y\right)^{\ast}$
is \textit{not always} uniformly elliptic in $y\in\mathbb{R}^{d}$
and $z\in\mathbb{R}_{0}^{m}$. Although Assumption \ref{Uniform elliptic}
is not satisfied, this case can be also discussed in our position
by ignoring any jump effects. Then, the sensitivity formulae are given
as follows (cf. \cite{Cass and Friz}): \begin{gather*}
\Gamma_{T}^{\left(1\right)}=\frac{L_{0,T}}{T},\quad\Gamma_{T}^{\left(2\right)}=L_{T}^{\varepsilon}+R_{T}^{\varepsilon}-Q_{T}^{\varepsilon},\quad\Gamma_{T}^{\left(3\right)}=\frac{L_{\tilde{T},T}^{\ast}\, L_{0,\tilde{T}}}{\tilde{T}^{2}}+\frac{1}{\tilde{T}}\int_{0}^{\tilde{T}}F^{1}\left(\varepsilon,s\right)dW_{s}.\end{gather*}
Moreover, remark that Assumption \ref{Levy measure} on the measure
$d\nu$ is not necessary. In case of $b_{z}\left(\varepsilon,y\right)\equiv0$,
these are exactly the Bismut-Elworthy-Li formulae. See \cite{Bismut-2}
and \cite{Elworthy and Li}. \hfill{}$\square$

}}\end{rem}

\begin{rem}\label{uniformly elliptic jump}{\rm{ Consider the case
where the $\mathbb{R}^{d}\otimes\mathbb{R}^{d}$-valued function $\partial_{z}b_{z}\left(\varepsilon,y\right)\partial_{z}b_{z}\left(\varepsilon,y\right)^{\ast}$
is uniformly elliptic in $y\in\mathbb{R}^{d}$ and $z\in\mathbb{R}_{0}^{m}$,
while the function $a\left(\varepsilon,y\right)a\left(\varepsilon,y\right)^{\ast}$
is \textit{not always} uniformly elliptic in $y\in\mathbb{R}^{d}$.
Although Assumption \ref{Uniform elliptic} is not satisfied, this
case can be also discussed in our position by ignoring any diffusion
effects. Then, the process can be of pure-jump type and of infinite-activity
type, and Assumption \ref{Levy measure} on $d\nu$ is essential.
The sensitivity formulae are\begin{align*}
\Gamma_{T}^{\left(1\right)} & =-\frac{V_{0,T}}{A_{0,T}}+\frac{K_{0,T}}{\left(A_{0,T}\right)^{2}},\quad\Gamma_{T}^{\left(2\right)}=-J_{T}^{\varepsilon},\\
\Gamma_{T}^{\left(3\right)} & =\left\{ -\frac{V_{\tilde{T},T}^{\ast}}{A_{\tilde{T},T}}+\frac{K_{\tilde{T},T}^{\ast}}{\left(A_{\tilde{T},T}\right)^{2}}+\frac{K_{\tilde{T},T}^{\ast}}{A_{\tilde{T},T}\, A_{0,\tilde{T}}}\right\} \left(-\frac{V_{0,\tilde{T}}}{A_{0,\tilde{T}}}+\frac{K_{0,\tilde{T}}}{\left(A_{0,\tilde{T}}\right)^{2}}\right)+\frac{K_{\tilde{T},T}^{\ast}\, K_{0,\tilde{T}}}{A_{\tilde{T},T}\,\left(A_{0,\tilde{T}}\ \right)^{3}}\\
 & \quad-\sum_{\sigma=2}^{3}\int_{0}^{\tilde{T}}\!\!\int_{\mathbb{R}_{0}^{m}}\frac{\mbox{div}_{z}\left[F^{\sigma}\left(\varepsilon,s,z\right)\right]\ }{A_{0,\tilde{T}}+\left|z\right|^{2}}d\hat{\mu}+\frac{1}{A_{0,\tilde{T}}}\sum_{\sigma=2}^{3}\int_{0}^{\tilde{T}}\!\!\int_{\mathbb{R}_{0}^{m}}\mbox{div}_{z}\left[F^{\sigma}\left(\varepsilon,s,z\right)\right]d\mu,\end{align*}
where $A_{\tau,t}=\int_{\tau}^{t}\!\!\int_{\mathbb{R}_{0}^{m}}\left|z\right|^{2}d\mu$.
\hfill{}$\square$

}}\end{rem}

\begin{rem}{\rm{ Bismut \cite{Bismut-1} obtained the integration
by parts formula for jump processes via the Girsanov transform, and
studied the existence of smooth densities. Then, it is crucial to
study the invertibility on the non-negative definite, symmetric matrices
valued random variable, which is called the Malliavin covariance matrix.
\textit{The H\"ormander type condition} on the coefficients of the
equation \eqref{eq:SDE}, instead of Assumption \ref{Uniform elliptic},
enables us to check the invertibility of the Malliavin covariance
matrix (cf. \cite{Komatsu-Takeuchi-1} and \cite{Komatsu-Takeuchi-2}).
Then, it would be possible to compute the concrete representations
as stated in Corollary \ref{sensitivity formulae} in \textit{the
hypoelliptic situation} via a similar manner to the one in the uniformly
elliptic situation, which will be studied elsewhere (cf. \cite{Takeuchi}).
\hfill{}$\square$

}}\end{rem}

\section{Proofs \label{sec:Proofs}}

We shall devote this section to prove our main results. For $t\in\left[0,T\right]$
and $\varphi\in C_{K}^{2}\left(\mathbb{R}^{d}\,;\,\mathbb{R}\right)$,
define \[
u\left(t,x\right)=\mathbb{E}\left[\varphi\left(x_{T-t}\right)|x_{0}=x\right].\]
Then, it holds that $u\in C_{b}^{1,2}\left([0,T)\times\mathbb{R}^{d}\,;\,\mathbb{R}\right)$,
$\lim_{t\nearrow T}u\left(t,x\right)=\varphi\left(x\right)$, and
$\left(\partial_{t}+\mathcal{L}^{\varepsilon}\right)u=0$ (cf. \cite{Gihman-Skorohod}).
The following lemma can be regarded as the martingale representation
on $\varphi\left(x_{T}\right)$, and plays a crucial role in what
follows. 

\begin{lem}\label{C-O}For $\varphi\in C_{K}^{2}\left(\mathbb{R}^{d}\,;\,\mathbb{R}\right)$,
it holds that \begin{equation}
\varphi\left(x_{T}\right)=\mathbb{E}\left[\varphi\left(x_{T}\right)\right]+\int_{0}^{T}\nabla u\left(s,x_{s}\right)a\left(\varepsilon,x_{s}\right)dW_{s}+\int_{0}^{T}\!\!\int_{\mathbb{R}_{0}^{m}}\mathfrak{B}_{z}^{\varepsilon}u\left(s,x_{s-}\right)d\tilde{\mu}.\label{eq:Ocone-Clark}\end{equation}
\end{lem}

\emph{Proof. }Let $t\in[0,T).$ Since $u\in C_{b}^{1,2}\left([0,T)\times\mathbb{R}^{d}\,;\,\mathbb{R}\right)$
and $\left(\partial_{t}+\mathcal{L}^{\varepsilon}\right)u=0$, the
It\^o formula yields\begin{equation}
u\left(t,x_{t}\right)=u\left(0,x\right)+\int_{0}^{t}\nabla u\left(s,x_{s}\right)a\left(\varepsilon,x_{s}\right)dW_{s}+\int_{0}^{t}\!\!\int_{\mathbb{R}_{0}^{m}}\mathfrak{B}_{z}^{\varepsilon}u\left(s,x_{s-}\right)d\tilde{\mu}.\label{eq:Clark-Ocone u}\end{equation}
Since $\varphi\in C_{K}^{2}\left(\mathbb{R}^{d}\,;\,\mathbb{R}\right)$,
it holds that $u\left(t,x_{t}\right)=\mathbb{E}\left[\varphi\left(x_{T}\right)|\mathcal{F}_{t}\right]\to\mathbb{E}\left[\varphi\left(x_{T}\right)|\mathcal{F}_{T}\right]=\varphi\left(x_{T}\right)$
as $t\nearrow T$ (cf. \cite{Ikeda-Watanabe}). It can be easily checked
from Assumption \ref{assumption: coefficients} that stochastic integrals
in the right hand side of the equality \eqref{eq:Clark-Ocone u} converge
to the ones in \eqref{eq:Ocone-Clark} as $t\nearrow T$, respectively.
\hfill{}$\square$

\bigskip{}

Taking the differential of \eqref{eq:Ocone-Clark} in Lemma \ref{C-O},
we have 

\begin{lem}\label{C-O-derivative}For $1\le k\le d$ and $\varphi\in C_{K}^{2}\left(\mathbb{R}^{d}\,;\,\mathbb{R}\right)$,
it holds that \begin{equation}
\begin{split}\nabla_{x_{k}}\left(\varphi\left(x_{T}\right)\right) & =\mathbb{E}\left[\nabla_{x_{k}}\left(\varphi\left(x_{T}\right)\right)\right]+\int_{0}^{T}\nabla_{x_{k}}\left(\nabla u\left(s,x_{s}\right)a\left(\varepsilon,x_{s}\right)\right)dW_{s}\\
 & \qquad+\int_{0}^{T}\!\!\int_{\mathbb{R}_{0}^{m}}\nabla_{x_{k}}\left(\mathfrak{B}_{z}^{\varepsilon}u\left(s,x_{s-}\right)\right)d\tilde{\mu}.\end{split}
\label{eq:1st derivative}\end{equation}
\end{lem}

\emph{Proof. }We shall write $x_{t}=x_{t}^{x}$, in order to emphasize
the dependence on $x\in\mathbb{R}^{d}$ throughout the proof. Taking
the derivative of the equality \eqref{eq:Ocone-Clark} in Lemma \ref{C-O},
we have\begin{align*}
\nabla_{x_{k}}\left(\varphi\left(x_{T}^{x}\right)\right) & =\nabla_{x_{k}}\left(\mathbb{E}\left[\varphi\left(x_{T}^{x}\right)\right]\right)+\nabla_{x_{k}}\left(\int_{0}^{T}\nabla u\left(s,x_{s}^{x}\right)a\left(\varepsilon,x_{s}^{x}\right)dW_{s}\right)\\
 & \qquad+\nabla_{x_{k}}\left(\int_{0}^{T}\!\!\int_{\mathbb{R}_{0}^{m}}\mathfrak{B}_{z}^{\varepsilon}u\left(s,x_{s-}^{x}\right)d\tilde{\mu}\right).\end{align*}
Let $0<\delta<1$, and $\bm{e}_{k}\in\mathbb{R}^{d}$ be the $k$-th
unit vector. Then, we have\[
\left|\frac{\mathbb{E}\left[\varphi\left(x_{T}^{x+\delta\bm{e}_{k}}\right)\right]-\mathbb{E}\left[\varphi\left(x_{T}^{x}\right)\right]}{\delta}-\mathbb{E}\left[\nabla_{x_{k}}\left(\varphi\left(x_{T}^{x}\right)\right)\right]\right|\le\int_{0}^{1}\mathbb{E}\left[\left|\nabla_{x+\sigma\delta\bm{e}_{k}}\left(\varphi\left(x_{T}^{x+\sigma\delta\bm{e}_{k}}\right)\right)-\nabla_{x_{k}}\left(\varphi\left(x_{T}^{x}\right)\right)\right|\right]d\sigma,\]
 which tends to $0$ as $\delta\searrow0$, because $\varphi\in C_{K}^{2}\left(\mathbb{R}^{d}\,;\,\mathbb{R}\right)$
and $Z_{t}\in\mathbb{L}^{2}\left(\Omega,\mathbb{P}\right)$. Hence,
we get \[
\nabla_{x_{k}}\left(\mathbb{E}\left[\varphi\left(x_{T}^{x}\right)\right]\right)=\mathbb{E}\left[\nabla_{x_{k}}\left(\varphi\left(x_{T}^{x}\right)\right)\right].\]
On the other hand, since $u\in C_{b}^{1,2}\left([0,T)\times\mathbb{R}^{d}\,;\,\mathbb{R}\right)$,
Assumption \ref{assumption: coefficients} and $x_{t},\, Z_{t}\in\mathbb{L}^{p}\left(\Omega,\mathbb{P}\right)$
for any $p>1$, we have \begin{align*}
 & \mathbb{E}\left[\left|\int_{0}^{T}\left(\frac{\nabla u\left(s,x_{s}^{x+\delta\bm{e}_{k}}\right)a\left(\varepsilon,x_{s}^{x+\delta\bm{e}_{k}}\right)-\nabla u\left(s,x_{s}^{x}\right)a\left(\varepsilon,x_{s}^{x}\right)}{\delta}-\nabla_{x_{k}}\left(\nabla u\left(s,x_{s}\right)a\left(\varepsilon,x_{s}\right)\right)\right)dW_{s}\right|^{2}\right]\\
 & \le\int_{0}^{T}\mathbb{E}\left[\int_{0}^{1}\left|\nabla_{x+\sigma\delta\bm{e}_{k}}\left(\nabla u\left(s,x_{s}^{x+\sigma\delta\bm{e}_{k}}\right)a\left(\varepsilon,x_{s}^{x+\sigma\delta\bm{e}_{k}}\right)\right)-\nabla_{x_{k}}\left(\nabla u\left(s,x_{s}^{x}\right)a\left(\varepsilon,x_{s}^{x}\right)\right)\right|^{2}d\sigma\right]ds,\end{align*}
which tends to $0$ as $\delta\searrow0$. Thus, we get \[
\nabla_{x_{k}}\left(\int_{0}^{T}\nabla u\left(s,x_{s}^{x}\right)a\left(\varepsilon,x_{s}^{x}\right)dW_{s}\right)=\int_{0}^{T}\nabla_{x_{k}}\left(\nabla u\left(s,x_{s}^{x}\right)a\left(\varepsilon,x_{s}^{x}\right)\right)dW_{s}.\]
Similarly to the above, it holds that \[
\nabla_{x_{k}}\left(\int_{0}^{T}\!\!\int_{\mathbb{R}_{0}^{m}}\mathfrak{B}_{z}^{\varepsilon}u\left(s,x_{s-}^{x}\right)d\tilde{\mu}\right)=\int_{0}^{T}\!\!\int_{\mathbb{R}_{0}^{m}}\nabla_{x_{k}}\left(\mathfrak{B}_{z}^{\varepsilon}u\left(s,x_{s-}^{x}\right)\right)d\tilde{\mu}.\]
\hfill{}$\square$

\begin{lem}\label{key eq} For $\varphi\in C_{K}^{2}\left(\mathbb{R}^{d}\,;\,\mathbb{R}\right)$,
it holds that\begin{gather*}
\mathbb{E}\left[\varphi\left(x_{T}\right)\,\int_{0}^{T}\!\!\int_{\mathbb{R}_{0}^{m}}\left|z\right|^{2}d\mu\right]=\mathbb{E}\left[\int_{0}^{T}\!\!\int_{\mathbb{R}_{0}^{m}}u\left(t,x_{t}+b_{z}\left(\varepsilon,x_{t}\right)\right)\left|z\right|^{2}d\hat{\mu}\right],\\
\mathbb{E}\left[\varphi\left(x_{T}\right)\,\int_{0}^{T}\!\!\int_{\mathbb{R}_{0}^{m}}h_{z}^{\varepsilon}\left(t\right)d\mu\right]=\mathbb{E}\left[\int_{0}^{T}\!\!\int_{\mathbb{R}_{0}^{m}}u\left(t,x_{t}+b_{z}\left(\varepsilon,x_{t}\right)\right)\, h_{z}^{\varepsilon}\left(t\right)d\hat{\mu}\right].\end{gather*}
 \end{lem}

\emph{Proof}. We shall only prove the second assertion, because the
first assertion can be proved in a similar manner. Write $M_{T}^{\varepsilon}=\int_{0}^{T}\!\!\int_{\mathbb{R}_{0}^{m}}h_{z}^{\varepsilon}\left(t\right)d\mu,\ \hat{M}_{T}^{\varepsilon}=\int_{0}^{T}\!\!\int_{\mathbb{R}_{0}^{m}}h_{z}^{\varepsilon}\left(t\right)\, d\hat{\mu}$.
Remark that \begin{align*}
 & \mathbb{E}\left[\int_{0}^{T}\nabla u\left(t,x_{t}\right)\, a\left(\varepsilon,x_{t}\right)\, dW_{t}\, M_{T}^{\varepsilon}\right]\\
 & =\mathbb{E}\left[\int_{0}^{T}\sum_{i=1}^{m}\nabla u\left(t,x_{t}\right)\, a_{i}\left(\varepsilon,x_{t}\right)\, M_{t}^{\varepsilon}\, dW_{t}^{i}\right]+\mathbb{E}\left[\int_{0}^{T}\left(\int_{0}^{t}\nabla u\left(s,x_{s}\right)\, a\left(\varepsilon,x_{s}\right)\, dW_{s}\right)\, dM_{t}^{\varepsilon}\right]\\
 & =\mathbb{E}\left[\int_{0}^{T}\left(\int_{0}^{t}\nabla u\left(s,x_{s}\right)\, a\left(\varepsilon,x_{s}\right)\, dW_{s}\right)\, d\hat{M}_{t}^{\varepsilon}\right],\end{align*}
and\begin{align*}
 & \mathbb{E}\left[\int_{0}^{T}\!\!\int_{\mathbb{R}_{0}^{m}}\mathfrak{B}_{z}^{\varepsilon}u\left(t,x_{t-}\right)\, d\tilde{\mu}\, M_{T}^{\varepsilon}\right]\\
 & =\mathbb{E}\left[\int_{0}^{T}\!\!\int_{\mathbb{R}_{0}^{m}}\mathfrak{B}_{z}^{\varepsilon}u\left(t,x_{t-}\right)\, h_{z}^{\varepsilon}\left(t\right)\, d\mu\right]+\mathbb{E}\left[\int_{0}^{T}\left(\int_{0}^{t}\!\!\int_{\mathbb{R}_{0}^{m}}\mathfrak{B}_{\theta}^{\varepsilon}u\left(s,x_{s-}\right)\, d\tilde{\mu}\right)\, dM_{t}^{\varepsilon}\right]\\
 & \qquad+\mathbb{E}\left[\int_{0}^{T}\!\!\int_{\mathbb{R}_{0}^{m}}\mathfrak{B}_{z}^{\varepsilon}u\left(t,x_{t-}\right)\, M_{t-}^{\varepsilon}\, d\tilde{\mu}\right]\\
 & =\mathbb{E}\left[\int_{0}^{T}\!\!\int_{\mathbb{R}_{0}^{m}}\mathfrak{B}_{z}^{\varepsilon}u\left(t,x_{t}\right)\, h_{z}^{\varepsilon}\left(t\right)\, d\hat{\mu}\right]+\mathbb{E}\left[\int_{0}^{T}\left(\int_{0}^{t}\!\!\int_{\mathbb{R}_{0}^{m}}\mathfrak{B}_{\theta}^{\varepsilon}u\left(s,x_{s-}\right)\, d\tilde{\mu}\right)\, d\hat{M}_{t}^{\varepsilon}\right]\end{align*}
from the It\^o formula. Since $u\in C_{b}^{1,2}\left([0,T)\times\mathbb{R}^{d}\,;\,\mathbb{R}\right)$,
$U_{t}\in\mathbb{L}^{p}\left(\Omega,\mathbb{P}\right)$ for any $p>1$,
and \begin{align*}
 & \mathbb{E}\left[\left\{ \int_{t}^{T}\nabla u\left(s,x_{s}\right)\, a\left(\varepsilon,x_{s}\right)\, dW_{s}+\int_{t}^{T}\!\!\int_{\mathbb{R}_{0}^{m}}\mathfrak{B}_{\theta}^{\varepsilon}u\left(s,x_{s-}\right)\, d\tilde{\mu}\right\} \, h_{z}^{\varepsilon}\left(t\right)\right]=\bm{0}\in\mathbb{R}^{l}\otimes\mathbb{R}^{d},\end{align*}
the equality \eqref{eq:Ocone-Clark} in Lemma \ref{C-O} enables us
to see that\begin{align*}
 & \mathbb{E}\left[\varphi\left(x_{T}\right)\, M_{T}^{\varepsilon}\right]\\
 & =\mathbb{E}\left[\varphi\left(x_{T}\right)\right]\,\mathbb{E}\left[\hat{M}_{T}^{\varepsilon}\right]+\mathbb{E}\left[\int_{0}^{T}\nabla u\left(t,x_{t}\right)\, a\left(\varepsilon,x_{t}\right)\, dW_{t}\, M_{T}^{\varepsilon}\right]+\mathbb{E}\left[\int_{0}^{T}\!\!\int_{\mathbb{R}_{0}^{m}}\mathfrak{B}_{z}^{\varepsilon}u\left(t,x_{t-}\right)\, d\tilde{\mu}\, M_{T}^{\varepsilon}\right]\\
 & =\mathbb{E}\left[\varphi\left(x_{T}\right)\,\hat{M}_{T}^{\varepsilon}\right]+\mathbb{E}\left[\int_{0}^{T}\!\!\int_{\mathbb{R}_{0}^{m}}\mathfrak{B}_{z}^{\varepsilon}u\left(t,x_{t}\right)\, h_{z}^{\varepsilon}\left(t\right)\, d\hat{\mu}\right]\\
 & \qquad-\mathbb{E}\left[\int_{0}^{T}\left(\int_{t}^{T}\nabla u\left(s,x_{s}\right)\, a\left(\varepsilon,x_{s}\right)\, dW_{s}\right)\, d\hat{M}_{t}^{\varepsilon}\right]-\mathbb{E}\left[\int_{0}^{T}\left(\int_{t}^{T}\!\!\int_{\mathbb{R}_{0}^{m}}\mathfrak{B}_{\theta}^{\varepsilon}u\left(s,x_{s-}\right)\, d\tilde{\mu}\right)\, d\hat{M}_{t}^{\varepsilon}\right]\\
 & =\mathbb{E}\left[\int_{0}^{T}\!\!\int_{\mathbb{R}_{0}^{m}}u^{\varepsilon}\left(t,x_{t}+b_{z}\left(\varepsilon,x_{t}\right)\right)\, h_{z}^{\varepsilon}\left(t\right)\, d\hat{\mu}\right].\end{align*}
\hfill{}$\square$

\subsection{Proofs of Theorem \ref{Theorem 1} and \ref{vega}\label{sub:Proof-of-Theorem1}}

We shall reveal each terms in Theorem \ref{Theorem 1} and \ref{vega}
in what follows. 

\begin{lem}\label{drift} Let $\varphi\in C_{K}^{2}\left(\mathbb{R}^{d}\,;\,\mathbb{R}\right)$.
Then, it holds that\begin{gather*}
\mathbb{E}\left[\nabla_{x}\left(\varphi\left(x_{T}\right)\right)\right]\, T=\mathbb{E}\left[\varphi\left(x_{T}\right)L_{0,T}\right],\\
\mathbb{E}\left[\nabla_{x}\left(\varphi\left(x_{T}\right)\right)\,\int_{0}^{T}f_{0}^{\varepsilon}\left(t\right)\, dt\right]=\mathbb{E}\left[\varphi\left(x_{T}\right)\, L_{T}^{\varepsilon}\right].\end{gather*}

\end{lem} 

\textit{Proof}. We shall only prove the second assertion, because
the first one can be done in a similar manner. Since Lemma \ref{C-O}
tells us that the process $\left\{ u\left(t,x_{t}\right)\,;\, t\in[0,T)\right\} $
is $\left(\mathcal{F}_{t}\right)$-martingale, so is $\left\{ \nabla_{x}\left(u\left(t,x_{t}\right)\right)\,;\, t\in[0,T)\right\} $,
similarly to the proof of Lemma \ref{C-O-derivative}. Then, for $t<\tau<T$,
we have \[
\mathbb{E}\left[\nabla_{x}\left(u\left(t,x_{t}\right)\right)f_{0}^{\varepsilon}\left(t\right)\right]=\mathbb{E}\left[\nabla_{x}\left(u\left(\tau,x_{\tau}\right)\right)f_{0}^{\varepsilon}\left(t\right)\right].\]
Hence, taking the limit as $\tau\nearrow T$ yields that $\mathbb{E}\left[\nabla_{x}\left(u\left(t,x_{t}\right)\right)f_{0}^{\varepsilon}\left(t\right)\right]=\mathbb{E}\left[\nabla_{x}\left(\varphi\left(x_{T}\right)\right)f_{0}^{\varepsilon}\left(t\right)\right]$,
because \begin{align*}
\nabla_{x}\left(u\left(\tau,x_{\tau}\right)\right) & =\nabla_{x}\left(u\left(0,x\right)\right)+\int_{0}^{\tau}\nabla_{x}\left(\nabla u\left(s,x_{s}\right)a\left(\varepsilon,x_{s}\right)\right)dW_{s}+\int_{0}^{\tau}\!\!\int_{\mathbb{R}_{0}^{m}}\nabla_{x}\left(\mathfrak{B}_{z}^{\varepsilon}u\left(s,x_{s-}\right)\right)d\tilde{\mu}\\
 & \to\nabla_{x}\left(u\left(0,x\right)\right)+\int_{0}^{T}\nabla_{x}\left(\nabla u\left(s,x_{s}\right)a\left(\varepsilon,x_{s}\right)\right)dW_{s}+\int_{0}^{T}\!\!\int_{\mathbb{R}_{0}^{m}}\nabla_{x}\left(\mathfrak{B}_{z}^{\varepsilon}u\left(s,x_{s-}\right)\right)d\tilde{\mu}\\
 & =\nabla_{x}\left(\varphi\left(x_{T}\right)\right).\end{align*}
 Therefore, the Fubini theorem and Lemma \ref{C-O} yield that \begin{align*}
\mathbb{E}\left[\nabla_{x}\left(\varphi\left(x_{T}\right)\right)\,\int_{0}^{T}f_{0}^{\varepsilon}\left(t\right)\, dt\right]\, & =\int_{0}^{T}\mathbb{E}\left[\nabla_{x}\left(u\left(t,x_{t}\right)\right)f_{0}^{\varepsilon}\left(t\right)\right]\, dt\\
 & =\mathbb{E}\left[\int_{0}^{T}\nabla u\left(t,x_{t}\right)\, a\left(\varepsilon,x_{t}\right)\, dW_{t}\, L_{T}^{\varepsilon}\right]=\mathbb{E}\left[\varphi\left(x_{T}\right)\, L_{T}^{\varepsilon}\right],\end{align*}
because of Assumption \ref{Uniform elliptic} and $U_{t}\in\mathbb{L}^{2}\left(\Omega,\mathbb{P}\right)$.
\hfill{}$\square$

\bigskip{}

The following lemma is the application of the integration by parts
formula on the Wiener space. 

\begin{lem}\label{diffusion} Let $\varphi\in C_{K}^{2}\left(\mathbb{R}^{d}\,;\,\mathbb{R}\right)$.
Then, it holds that\[
\mathbb{E}\left[\nabla_{x}\left(\varphi\left(x_{T}\right)\right)\,\int_{0}^{T}f^{\varepsilon}\left(t\right)\circ dW_{t}\right]=\mathbb{E}\left[\varphi\left(x_{T}\right)\,\left(R_{T}^{\varepsilon}-Q_{T}^{\varepsilon}\right)\right].\]
 \end{lem} 

\textit{Proof}. Since $D_{s}\varphi\left(x_{T}\right)=\nabla\varphi\left(x_{T}\right)\, Z_{T}\, U_{s}\, a\left(\varepsilon,x_{s}\right)$
for $s\in\left[0,T\right]$ from the chain rule on the operator $D$,
the integration by parts formula implies that\begin{align*}
\mathbb{E}\left[\nabla_{x}\left(\varphi\left(x_{T}\right)\right)\, G_{T}\right] & =\mathbb{E}\left[\frac{1}{T}\int_{0}^{T}D_{s}\varphi\left(x_{T}\right)\, a\left(\varepsilon,x_{s}\right)^{-1}\, Z_{s}\, G_{T}\, ds\right]\\
 & =\mathbb{E}\left[\varphi\left(x_{T}\right)\ \frac{1}{T}\, D^{\ast}\left(a\left(\varepsilon,x_{\cdot}\right)^{-1}\, Z_{\cdot}\, G_{T}\right)\right],\end{align*}
where $D^{\ast}$ is the Skorokhod integral operator. Remark that
$G_{T}\in\mathbb{D}_{\infty}\left(\mathbb{R}^{l}\otimes\mathbb{R}^{d}\right)$
from Assumption \ref{assumption: coefficients} (cf. \cite{Nualart}).
Then, we see that\begin{align*}
D^{\ast}\left(a\left(\varepsilon,x_{\cdot}\right)^{-1}\, Z_{\cdot}\, G_{T}\right) & =D^{\ast}\left(a\left(\varepsilon,x_{\cdot}\right)^{-1}\, Z_{\cdot}\right)\, G_{T}-\int_{0}^{T}\mbox{Tr}\left[a\left(\varepsilon,x_{s}\right)^{-1}\, Z_{s}\, D_{s}G_{T}\right]\, ds\\
 & =\left\{ \int_{0}^{T}\left(dW_{s}\right)^{\ast}\, a\left(\varepsilon,x_{s}\right)^{-1}\, Z_{s}\right\} \, G_{T}-\int_{0}^{T}\mbox{Tr}\left[a\left(\varepsilon,x_{s}\right)^{-1}\, Z_{s}\, D_{s}G_{T}\right]\, ds\\
 & =T\, R_{T}^{\varepsilon}-T\, Q_{T}^{\varepsilon}\end{align*}
from Proposition I-1.3.3 in \cite{Nualart}. \hfill{}$\square$

\begin{lem}\label{discontinuous} Let $\varphi\in C_{K}^{2}\left(\mathbb{R}^{d}\,;\,\mathbb{R}\right)$.
Then, it holds that\begin{gather*}
\mathbb{E}\left[\nabla_{x}\left(\varphi\left(x_{T}\right)\right)\,\int_{0}^{T}\!\!\int_{\mathbb{R}_{0}^{m}}\left|z\right|^{2}d\mu\right]=-\mathbb{E}\left[\varphi\left(x_{T}\right)J_{0,T}\right],\\
\mathbb{E}\left[\nabla_{x}\left(\varphi\left(x_{T}\right)\right)\,\int_{0}^{T}\!\!\int_{\mathbb{R}_{0}^{m}}h_{z}^{\varepsilon}\left(t\right)\, d\mu\right]=-\mathbb{E}\left[\varphi\left(x_{T}\right)\, J_{T}^{\varepsilon}\right].\end{gather*}
 \end{lem} 

\textit{Proof}. We shall prove the second assertion only, because
the first assertion can be obtained in a similar manner. Recall $M_{T}^{\varepsilon}=\int_{0}^{T}\!\!\int_{\mathbb{R}_{0}^{m}}h_{z}^{\varepsilon}\left(t\right)d\mu$
and $\hat{M}_{T}^{\varepsilon}=\int_{0}^{T}\!\!\int_{\mathbb{R}_{0}^{m}}h_{z}^{\varepsilon}\left(t\right)\, d\hat{\mu}$.
Lemma \ref{C-O} implies that\begin{align*}
\mathbb{E}\left[\int_{0}^{T}\!\!\int_{\mathbb{R}_{0}^{m}}\mathfrak{B}_{z}^{\varepsilon}u\left(t,x_{t}\right)\,\nabla_{x}\left(h_{z}^{\varepsilon}\left(t\right)\right)\, d\hat{\mu}\right] & =\mathbb{E}\left[\int_{0}^{T}\!\!\int_{\mathbb{R}_{0}^{m}}\mathfrak{B}_{z}^{\varepsilon}u\left(t,x_{t-}\right)\, d\tilde{\mu}\,\int_{0}^{T}\!\!\int_{\mathbb{R}_{0}^{m}}\nabla_{x}\left(h_{z}^{\varepsilon}\left(t\right)\right)\, d\tilde{\mu}\right]\\
 & =\mathbb{E}\left[\varphi\left(x_{T}\right)\,\int_{0}^{T}\!\!\int_{\mathbb{R}_{0}^{m}}\nabla_{x}\left(h_{z}^{\varepsilon}\left(t\right)\right)\, d\tilde{\mu}\right].\end{align*}
On the other hand, it holds that \[
\mathbb{E}\left[\int_{0}^{T}\!\!\int_{\mathbb{R}_{0}^{m}}u\left(t,x_{t}\right)\,\nabla_{x}\left(h_{z}^{\varepsilon}\left(t\right)\right)\, d\hat{\mu}\right]=\mathbb{E}\left[\varphi\left(x_{T}\right)\int_{0}^{T}\!\!\int_{\mathbb{R}_{0}^{m}}\nabla_{x}\left(h_{z}^{\varepsilon}\left(t\right)\right)\, d\hat{\mu}\right]\]
from the Fubini theorem. Thus, we have\begin{align*}
\mathbb{E}\left[\int_{0}^{T}\!\!\int_{\mathbb{R}_{0}^{m}}u\left(t,x_{t}+b_{z}\left(\varepsilon,x_{t}\right)\right)\,\nabla_{x}\left(h_{z}^{\varepsilon}\left(t\right)\right)\, d\hat{\mu}\right] & =\mathbb{E}\left[\varphi\left(x_{T}\right)\int_{0}^{T}\!\!\int_{\mathbb{R}_{0}^{m}}\nabla_{x}\left(h_{z}^{\varepsilon}\left(t\right)\right)\, d\mu\right]\\
 & =\mathbb{E}\left[\varphi\left(x_{T}\right)\,\nabla_{x}M_{T}^{\varepsilon}\right].\end{align*}
Here, the second equality can be justified, similarly to the proof
of Lemma \ref{C-O-derivative}. Furthermore, multiplying $J_{T}^{\varepsilon}$
by the equality \eqref{eq:Ocone-Clark} in Lemma \ref{C-O}, it holds
that \begin{align*}
\mathbb{E}\left[\varphi\left(x_{T}\right)\, J_{T}^{\varepsilon}\right] & =\mathbb{E}\left[\int_{0}^{T}\!\!\int_{\mathbb{R}_{0}^{m}}\mathfrak{B}_{z}^{\varepsilon}u\left(s,x_{s}\right)\,\mbox{{\rm {div}}}_{z}\left[g\left(z\right)\,\tilde{v}_{s}\left(\varepsilon,z\right)\right]\, dz\, ds\right]\\
 & =-\mathbb{E}\left[\int_{0}^{T}\!\!\int_{\mathbb{R}_{0}^{m}}\partial_{z}\left(u\left(s,x_{s}+b_{z}\left(\varepsilon,x_{s}\right)\right)\right)\,\tilde{v}_{s}\left(\varepsilon,z\right)\, d\hat{\mu}\right]\\
 & =-\mathbb{E}\left[\int_{0}^{T}\!\!\int_{\mathbb{R}_{0}^{m}}\nabla_{x}\left(u\left(s,x_{s}+b_{z}\left(\varepsilon,x_{s}\right)\right)\right)\, h_{z}^{\varepsilon}\left(s\right)\, d\hat{\mu}\right].\end{align*}
Here we have used the integration by parts for the second equality,
from the condition (iii) in Assumption \ref{Levy measure}. Therefore,
Assumption \ref{assumption: coefficients} and $\varphi\in C_{K}^{2}\left(\mathbb{R}^{d}\,;\,\mathbb{R}\right)$
enables us to obtain that \begin{align*}
\mathbb{E}\left[\nabla_{x}\left(\varphi\left(x_{T}\right)\right)\, M_{T}^{\varepsilon}\right] & =\mathbb{E}\left[\nabla_{x}\left(\varphi\left(x_{T}\right)\, M_{T}^{\varepsilon}\right)\right]-\mathbb{E}\left[\varphi\left(x_{T}\right)\,\nabla_{x}M_{T}^{\varepsilon}\right]\\
 & =\nabla_{x}\left(\mathbb{E}\left[\int_{0}^{T}\!\!\int_{\mathbb{R}_{0}^{m}}u\left(t,x_{t}+b_{z}\left(\varepsilon,x_{t}\right)\right)\, h_{z}^{\varepsilon}\left(t\right)\, d\hat{\mu}\right]\right)\\
 & \qquad-\mathbb{E}\left[\int_{0}^{T}\!\!\int_{\mathbb{R}_{0}^{m}}u\left(t,x_{t}+b_{z}\left(\varepsilon,x_{t}\right)\right)\,\nabla_{x}\left(h_{z}^{\varepsilon}\left(t\right)\right)\, d\hat{\mu}\right]\\
 & =\mathbb{E}\left[\int_{0}^{T}\!\!\int_{\mathbb{R}_{0}^{m}}\nabla_{x}\left(u\left(s,x_{s}+b_{z}\left(\varepsilon,x_{s}\right)\right)\right)\, h_{z}^{\varepsilon}\left(s\right)\, d\hat{\mu}\right]\\
 & =-\mathbb{E}\left[\varphi\left(x_{T}\right)\, J_{T}^{\varepsilon}\right].\end{align*}
\hfill{}$\square$

\begin{cor}\label{Corollary 1} For $\varphi\in C_{K}^{2}\left(\mathbb{R}^{d}\,;\,\mathbb{R}\right)$,
it holds that

\begin{equation}
\mathbb{E}\left[\nabla_{x}\left(\varphi\left(x_{T}\right)\right)A_{0,T}\right]=\mathbb{E}\left[\varphi\left(x_{T}\right)\left(L_{0,T}-J_{0,T}\right)\right].\label{eq:equality}\end{equation}
\end{cor}\bigskip{}

\emph{Proof of Theorem \ref{Theorem 1}.} Our goal is to get rid of
$A_{0,T}$ from the left hand side of the equality in Corollary \ref{Corollary 1}.
The Fubini theorem leads to\[
\mathbb{E}\left[\varphi\left(x_{T}\right)\right]=\mathbb{E}\left[\frac{\varphi\left(x_{T}\right)A_{0,T}}{\ A_{0,T}\ }\right]=\int_{0}^{\infty}\mathbb{E}\left[\varphi\left(x_{T}\right)A_{0,T}\,\exp\left(-\lambda\, A_{0,T}-N_{0,T}^{\lambda}\right)\right]e^{N_{0,T}^{\lambda}}d\lambda,\]
where $N_{\tau,t}^{\lambda}=\int_{\tau}^{t}\!\!\int_{\mathbb{R}_{0}^{m}}\left(e^{-\lambda\left|z\right|^{2}}-1\right)\, d\hat{\mu}$.
Define a new probability measure $\mathbb{P}^{\lambda}$ via the Girsanov
transform \[
\frac{d\mathbb{P}^{\lambda}}{d\mathbb{P}}\Big|_{\mathcal{F}_{T}}=\exp\left\{ -\int_{0}^{T}\!\!\int_{\mathbb{R}_{0}^{m}}\lambda\left|z\right|^{2}d\mu-N_{0,T}^{\lambda}\right\} ,\]
and denote by $\mathbb{E}^{\lambda}\left[\ \cdot\ \right]$ the expectation
with respect to the measure $\mathbb{P}^{\lambda}$. Then, under the
measure $\mathbb{P}^{\lambda}$, $d\mu$ is the Poisson random measure
with the intensity $d\hat{\mu}_{\lambda}:=\exp\left[-\lambda\left|z\right|^{2}\right]d\hat{\mu}$,
and $d\tilde{\mu}_{\lambda}:=d\mu-d\hat{\mu}_{\lambda}$ is a martingale
measure. See \cite{Kunita} for details. 

We shall rewrite the equation \eqref{eq:SDE} as follows: \[
dx_{t}=\tilde{a}_{0}\left(\varepsilon,x_{t}\right)\, dt+a\left(\varepsilon,x_{t}\right)\circ dW_{t}+\int_{\mathbb{R}_{0}^{m}}b_{z}\left(\varepsilon,x_{t-}\right)\, d\overline{\mu}_{\lambda},\]
where \[
\tilde{a}_{0}\left(\varepsilon,y\right)=a_{0}\left(\varepsilon,y\right)+\int_{\left|z\right|\le1}b_{z}\left(\varepsilon,y\right)\,\left(e^{-\lambda\left|z\right|^{2}}-1\right)\, d\nu,\quad d\overline{\mu}_{\lambda}=\mathbb{I}_{\left(\left|z\right|\le1\right)}\, d\tilde{\mu}_{\lambda}+\mathbb{I}_{\left(\left|z\right|>1\right)}\, d\mu.\]
In a similar manner to Corollary \ref{Corollary 1}, we can get \[
\nabla_{x}\left(\mathbb{E}^{\lambda}\left[\varphi\left(x_{T}\right)A_{0,T}\right]\right)=\mathbb{E}^{\lambda}\left[\varphi\left(x_{T}\right)\left(L_{0,T}-J_{0,T}^{\left(\lambda\right)}\right)\right],\]
where \[
J_{\tau,t}^{\left(\lambda\right)}=\int_{\tau}^{t}\!\!\int_{\mathbb{R}_{0}^{m}}\frac{\mbox{div}_{z}\left[e^{-\lambda\left|z\right|^{2}}g\left(z\right)v_{s}\left(\varepsilon,z\right)\right]}{e^{-\lambda\left|z\right|^{2}}g\left(z\right)}\, d\tilde{\mu}_{\lambda}\]
 for $0\le\tau\le t\le T$. The Fubini theorem yields that \[
\begin{split}\int_{0}^{\infty}e^{-\lambda T+N_{0,T}^{\lambda}}\,\mathbb{E}^{\lambda}\left[\varphi\left(x_{T}\right)L_{0,T}\right]d\lambda & =\mathbb{E}\left[\left\{ \int_{0}^{\infty}e^{-\lambda A_{0,T}}d\lambda\right\} \varphi\left(x_{T}\right)L_{0,T}\right]=\mathbb{E}\left[\varphi\left(x_{T}\right)\frac{L_{0,T}}{\ A_{0,T}\ \ }\right].\end{split}
\]
Since $d\tilde{\mu}_{\lambda}=d\tilde{\mu}+\left(1-e^{-\lambda\left|z\right|^{2}}\right)d\hat{\mu}$
and \begin{gather*}
\frac{\mbox{div}_{z}\left[e^{-\lambda\left|z\right|^{2}}g\left(z\right)v_{s}\left(\varepsilon,z\right)\right]}{e^{-\lambda\left|z\right|^{2}}g\left(z\right)}=-2\lambda z^{\ast}v_{s}\left(\varepsilon,z\right)+\frac{\ \mbox{div}_{z}\left[g\left(z\right)v_{s}\left(\varepsilon,z\right)\right]}{g\left(z\right)},\end{gather*}
we have $J_{\tau,t}^{\left(\lambda\right)}=J_{\tau,t}-\lambda K_{\tau,t}$
from the condition (iii) in Assumption \ref{Levy measure}. Hence,
we can get \[
\begin{split}\int_{0}^{\infty}e^{-\lambda T+N_{0,T}^{\lambda}}\,\mathbb{E}^{\lambda}\left[\varphi\left(x_{T}\right)J_{0,T}^{\left(\lambda\right)}\right]d\lambda & =\mathbb{E}\left[\int_{0}^{\infty}e^{-\lambda A_{0,T}}\varphi\left(x_{T}\right)\left(J_{0,T}-\lambda K_{0,T}\right)d\lambda\right]\\
 & =\mathbb{E}\left[\varphi\left(x_{T}\right)\frac{J_{0,T}}{A_{0,T}\ }\right]-\mathbb{E}\left[\varphi\left(x_{T}\right)\frac{K_{0,T}}{\left(A_{0,T}\right)^{2}}\right]\end{split}
\]
from the Fubini theorem. The proof of Theorem \ref{Theorem 1} is
complete. \hfill{}$\square$

\textit{\medskip{}
}

\textit{Proof of Theorem \ref{vega}.} By summing up the equalities
in Lemma \ref{drift}, \ref{diffusion} and \ref{discontinuous},
the assertion of Theorem \ref{vega} holds. \hfill{}$\square$

\subsection{Proof of Theorem \ref{Theorem 2}\label{sub:Proof-of-Theorem2}}

We shall reveal each terms in Theorem \ref{Theorem 2}. Concerning
the continuous part, it holds that 

\begin{lem}\label{Continuous 2nd}For $1\le j,\ k\le d$, and $\varphi\in C_{K}^{2}\left(\mathbb{R}^{d}\,;\,\mathbb{R}\right)$,
it holds that \begin{equation}
\tilde{T}\,\nabla_{x_{j}}\nabla_{x_{k}}\left(\mathbb{E}\left[\varphi\left(x_{T}\right)\right]\right)=\mathbb{E}\left[\varphi\left(x_{T}\right)\left\{ \Gamma_{\tilde{T},T}^{\left(1\right),j}\, L_{0,\tilde{T}}^{k}+\int_{0}^{\tilde{T}}\!\sum_{i=1}^{m}F_{ijk}^{1}\left(\varepsilon,s\right)dW_{s}^{i}\right\} \right].\label{eq:2nd deri conti}\end{equation}
\end{lem}

\emph{Proof.} Define $\ell_{s}\left(\varepsilon\right)=a\left(\varepsilon,x_{s}\right)^{-1}Z_{s}$.
Multiplying $L_{0,T}^{k}$ by \eqref{eq:1st derivative} in Lemma
\ref{C-O-derivative}, we have\begin{align*}
\mathbb{E}\left[\left(\nabla_{x_{j}}\left(\varphi\left(x_{T}\right)\right)\right)L_{0,T}^{k}\right] & =\mathbb{E}\left[\int_{0}^{T}\sum_{i=1}^{m}\nabla_{x_{j}}\left\{ \nabla u\left(s,x_{s}\right)a_{i}\left(\varepsilon,x_{s}\right)\right\} \,\ell_{s}^{ik}\left(\varepsilon\right)ds\right]\\
 & =\mathbb{E}\left[\int_{0}^{T}\sum_{\beta,\gamma=1}^{d}\nabla_{\gamma}\nabla_{\beta}u\left(s,x_{s}\right)Z_{s}^{\gamma j}Z_{s}^{\beta k}ds\right]\\
 & \qquad+\mathbb{E}\left[\int_{0}^{T}\sum_{i=1}^{m}\sum_{\beta=1}^{d}\nabla_{\beta}u\left(s,x_{s}\right)\nabla_{x_{j}}\left(a_{i}^{\beta}\left(\varepsilon,x_{s}\right)\right)\ell_{s}^{ik}\left(\varepsilon\right)ds\right]\\
 & =\mathbb{E}\left[\int_{0}^{T}\nabla_{x_{j}}\nabla_{x_{k}}\left(u\left(s,x_{s}\right)\right)ds\right]-\mathbb{E}\left[\int_{0}^{T}\sum_{\beta=1}^{d}\nabla_{\beta}u\left(s,x_{s}\right)\nabla_{x_{k}}Z_{s}^{\beta j}ds\right]\\
 & \qquad-\mathbb{E}\left[\int_{0}^{T}\sum_{i=1}^{m}\sum_{\beta=1}^{d}\nabla u\left(s,x_{s}\right)a_{i}\left(\varepsilon,x_{s}\right)\nabla_{x_{j}}\left[a\left(\varepsilon,x_{s}\right)^{-1}\right]_{i\beta}Z_{s}^{\beta k}ds\right],\end{align*}
because $\sum_{i=1}^{m}\nabla_{\gamma}a_{i}^{\beta}\left(y\right)\left[a\left(y\right)^{-1}\right]_{i\delta}=-\sum_{i=1}^{m}a_{i}^{\beta}\left(y\right)\nabla_{\gamma}\left[a\left(y\right)^{-1}\right]_{i\delta}$.
Similarly to the proof of Lemma \ref{C-O-derivative}, we can get
$\mathbb{E}\left[\nabla_{x_{j}}\nabla_{x_{k}}\left(u\left(s,x_{s}\right)\right)\right]=\mathbb{E}\left[\nabla_{x_{j}}\nabla_{x_{k}}\left(u\left(t,x_{t}\right)\right)\right]$
for $0\le s\le t\le T$. Since $\varphi\in C_{K}^{2}\left(\mathbb{R}^{d}\,;\,\mathbb{R}\right)$,
taking the limit as $t\nearrow T$ leads to \[
\mathbb{E}\left[\nabla_{x_{j}}\nabla_{x_{k}}\left(u\left(s,x_{s}\right)\right)\right]=\mathbb{E}\left[\nabla_{x_{j}}\nabla_{x_{k}}\left(\varphi\left(x_{T}\right)\right)\right]=\nabla_{x_{j}}\nabla_{x_{k}}\left(\mathbb{E}\left[\varphi\left(x_{T}\right)\right]\right).\]
Thus, we see that \begin{equation}
\begin{split}T\,\nabla_{x_{j}}\nabla_{x_{k}}\left(\mathbb{E}\left[\varphi\left(x_{T}\right)\right]\right) & =\mathbb{E}\left[\nabla_{x_{j}}\left(\varphi\left(x_{T}\right)\right)\, L_{0,T}^{k}\right]+\mathbb{E}\left[\int_{0}^{T}\sum_{\beta=1}^{d}\nabla_{\beta}u\left(s,x_{s}\right)\nabla_{x_{k}}Z_{s}^{\beta j}ds\right]\\
 & \quad+\mathbb{E}\left[\int_{0}^{T}\sum_{i=1}^{m}\sum_{\beta=1}^{d}\nabla u\left(s,x_{s}\right)a_{i}\left(\varepsilon,x_{s}\right)\nabla_{x_{j}}\left[a\left(\varepsilon,x_{s}\right)^{-1}\right]_{i\beta}Z_{s}^{\beta k}ds\right].\end{split}
\label{eq:2nd derivative(continuous part)}\end{equation}

Denote by $\left\{ P_{t}\,;\, t\in\left[0,T\right]\right\} $ the
$\left(C_{0}\right)$-semigroup associated with the process $\left\{ x_{t}\,;\, t\in\left[0,T\right]\right\} $.
We shall replace $T$ and $\varphi$ in the equality \eqref{eq:2nd derivative(continuous part)}
by $\tilde{T}$ and $P_{\tilde{T}}\varphi$, respectively. Then, it
holds that \[
\tilde{T}\,\nabla_{x_{j}}\nabla_{x_{k}}\left(\mathbb{E}\left[P_{\tilde{T}}\varphi\left(x_{\tilde{T}}\right)\right]\right)=\tilde{T}\,\nabla_{x_{j}}\nabla_{x_{k}}\left(\mathbb{E}\left[\mathbb{E}\left[\varphi\left(x_{T}\right)\big|\mathcal{F}_{\tilde{T}}\right]\right]\right)=\tilde{T}\,\nabla_{x_{j}}\nabla_{x_{k}}\left(\mathbb{E}\left[\varphi\left(x_{T}\right)\right]\right).\]
In a similar manner to Theorem \ref{Theorem 1}, the first term of
the right hand side of \eqref{eq:2nd derivative(continuous part)}
is equal to\[
\mathbb{E}\left[\nabla_{x_{j}}\left(P_{\tilde{T}}\varphi\left(x_{\tilde{T}}\right)\right)\, L_{0,\tilde{T}}^{k}\right]=\mathbb{E}\left[\nabla_{x_{j}}\left(\mathbb{E}\left[\varphi\left(x_{T}\right)\big|\mathcal{F}_{\tilde{T}}\right]\,\right)L_{0,\tilde{T}}^{k}\right]=\mathbb{E}\left[\varphi\left(x_{T}\right)\Gamma_{\tilde{T},T}^{\left(1\right),j}\, L_{0,\tilde{T}}^{k}\right].\]
 Define $\tilde{u}\left(t,x\right)=\mathbb{E}\left[P_{\tilde{T}}\varphi\left(x_{\tilde{T}}\right)\big|x_{0}=x\right]$
for $t\in[0,\tilde{T}]$ and $x\in\mathbb{R}^{d}$. Replacing $T$
and $\varphi$ in Lemma \ref{C-O} by $\tilde{T}$ and $P_{\tilde{T}}\varphi$,
respectively, we have \begin{equation}
\begin{split}P_{\tilde{T}}\varphi\left(x_{\tilde{T}}\right) & =\mathbb{E}\left[P_{\tilde{T}}\varphi\left(x_{\tilde{T}}\right)\right]+\int_{0}^{\tilde{T}}\nabla\tilde{u}\left(s,x_{s}\right)a\left(x_{s}\right)dW_{s}+\int_{0}^{\tilde{T}}\!\!\int_{\mathbb{R}_{0}^{m}}\mathfrak{B}_{z}\tilde{u}\left(s,x_{s-}\right)d\tilde{\mu}.\end{split}
\label{eq:semi-group}\end{equation}
 Multiplying $\int_{0}^{\tilde{T}}\sum_{i=1}^{m}\sum_{\beta=1}^{d}\left[a\left(\varepsilon,x_{s}\right)^{-1}\right]_{i\beta}\nabla_{x_{k}}Z_{s}^{\beta j}dW_{s}^{i}$
by the equality \eqref{eq:semi-group}, we have\begin{align*}
 & \mathbb{E}\left[\varphi\left(x_{T}\right)\int_{0}^{\tilde{T}}\sum_{i=1}^{m}\sum_{\beta=1}^{d}\left[a\left(\varepsilon,x_{s}\right)^{-1}\right]_{i\beta}\nabla_{x_{k}}Z_{s}^{\beta j}dW_{s}^{i}\right]\\
 & =\mathbb{E}\left[P_{\tilde{T}}\varphi\left(x_{\tilde{T}}\right)\int_{0}^{\tilde{T}}\sum_{i=1}^{m}\sum_{\beta=1}^{d}\left[a\left(\varepsilon,x_{s}\right)^{-1}\right]_{i\beta}\nabla_{x_{k}}Z_{s}^{\beta j}dW_{s}^{i}\right]\\
 & =\mathbb{E}\left[\int_{0}^{\tilde{T}}\sum_{\beta=1}^{d}\nabla_{\beta}\tilde{u}\left(s,x_{s}\right)\nabla_{x_{k}}Z_{s}^{\beta j}ds\right].\end{align*}
Multiplying $\int_{0}^{\tilde{T}}\sum_{i=1}^{m}\sum_{\beta=1}^{d}\nabla_{x_{j}}\left[a\left(\varepsilon,x_{s}\right)^{-1}\right]_{i\beta}Z_{s}^{\beta k}dW_{s}^{i}$
by the equality \eqref{eq:semi-group}, we have\[
\begin{split} & \mathbb{E}\left[\varphi\left(x_{T}\right)\int_{0}^{\tilde{T}}\sum_{i=1}^{m}\sum_{\beta=1}^{d}\nabla_{x_{j}}\left[a\left(\varepsilon,x_{s}\right)^{-1}\right]_{i\beta}\, Z_{s}^{\beta k}dW_{s}^{i}\right]\\
 & =\mathbb{E}\left[P_{\tilde{T}}\varphi\left(x_{\tilde{T}}\right)\int_{0}^{\tilde{T}}\sum_{i=1}^{m}\sum_{\beta=1}^{d}\nabla_{x_{j}}\left[a\left(\varepsilon,x_{s}\right)^{-1}\right]_{i\beta}\, Z_{s}^{\beta k}dW_{s}^{i}\right]\\
 & =\mathbb{E}\left[\int_{0}^{\tilde{T}}\sum_{i=1}^{m}\sum_{\beta=1}^{d}\nabla\tilde{u}\left(s,x_{s}\right)a_{i}\left(x_{s}\right)\nabla_{x_{j}}\left[a\left(\varepsilon,x_{s}\right)^{-1}\right]_{i\beta}Z_{s}^{\beta k}ds\right].\end{split}
\]
Therefore we can get\[
\begin{split}\tilde{T}\,\nabla_{x_{j}}\nabla_{x_{k}}\left(\mathbb{E}\left[\varphi\left(x_{T}\right)\right]\right) & =\mathbb{E}\left[\varphi\left(x_{T}\right)\Gamma_{\tilde{T},T}^{\left(1\right),j}L_{0,\tilde{T}}^{k}\right]+\mathbb{E}\left[\varphi\left(x_{T}\right)\int_{0}^{\tilde{T}}\sum_{i=1}^{m}\sum_{\beta=1}^{d}\left[a\left(\varepsilon,x_{s}\right)^{-1}\right]_{i\beta}\nabla_{x_{k}}Z_{s}^{\beta j}dW_{s}^{i}\right]\\
 & \qquad+\mathbb{E}\left[\varphi\left(x_{T}\right)\int_{0}^{\tilde{T}}\sum_{i=1}^{m}\sum_{\beta=1}^{d}\nabla_{x_{j}}\left[a\left(\varepsilon,x_{s}\right)^{-1}\right]_{i\beta}\, Z_{s}^{\beta k}dW_{s}^{i}\right]\\
 & =\mathbb{E}\left[\varphi\left(x_{T}\right)\left\{ \Gamma_{\tilde{T},T}^{\left(1\right),j}\, L_{0,\tilde{T}}^{k}+\int_{0}^{\tilde{T}}\sum_{i=1}^{m}F_{ijk}^{1}\left(\varepsilon,s\right)dW_{s}^{i}\right\} \right].\end{split}
\]
\hfill{}$\square$\bigskip{}

Concerning the jump part, it holds that 

\begin{lem}\label{Jump 2nd}For $1\le j,\ k\le d$, and $\varphi\in C_{K}^{2}\left(\mathbb{R}^{d}\,;\,\mathbb{R}\right)$,
it holds that\begin{equation}
\begin{aligned} & \nabla_{x_{j}}\nabla_{x_{k}}\left(\mathbb{E}\left[\varphi\left(x_{T}\right)\int_{0}^{\tilde{T}}\!\!\int_{\mathbb{R}_{0}^{m}}\left|z\right|^{2}d\mu\right]\right)\\
 & =\mathbb{E}\left[\varphi\left(x_{T}\right)\left\{ -\Gamma_{\tilde{T},T}^{\left(1\right),j}J_{0,\tilde{T}}^{k}+\int_{0}^{\tilde{T}}\!\!\int_{\mathbb{R}_{0}^{m}}\sum_{\sigma=2}^{3}\left\{ \mbox{div}_{z}\left[F^{\sigma}\left(\varepsilon,s,z\right)\right]\right\} _{jk}d\tilde{\mu}\right\} \right].\end{aligned}
\label{eq:2nd derivative(jump part)}\end{equation}
\end{lem}

\emph{Proof.} Since\begin{align*}
 & \nabla_{x_{j}}\nabla_{x_{k}}\left(u\left(s,x_{s}+b_{z}\left(\varepsilon,x_{s}\right)\right)\right)\\
 & =\nabla_{x_{k}}\left(\sum_{\beta=1}^{d}\nabla_{\beta}u\left(s,x_{s}+b_{z}\left(\varepsilon,x_{s}\right)\right)\left[\left(I_{d}+\nabla b_{z}\left(\varepsilon,x_{s}\right)\right)Z_{s}\right]_{\beta j}\right)\\
 & =\sum_{\beta,\delta=1}^{d}\nabla_{\delta}\nabla_{\beta}u\left(s,x_{s}+b_{z}\left(\varepsilon,x_{s}\right)\right)\left[\left(I_{d}+\nabla b_{z}\left(\varepsilon,x_{s}\right)\right)Z_{s}\right]_{\beta j}\left[\left(I_{d}+\nabla b_{z}\left(\varepsilon,x_{s}\right)\right)Z_{s}\right]_{\delta k}\\
 & \qquad+\sum_{\beta=1}^{d}\nabla_{\beta}u\left(s,x_{s}+b_{z}\left(\varepsilon,x_{s}\right)\right)\nabla_{x_{k}}\left(\left[\left(I_{d}+\nabla b_{z}\left(\varepsilon,x_{s}\right)\right)Z_{s}\right]_{\beta j}\right),\end{align*}
multiplying $J_{0,T}^{k}$ by the equality \eqref{eq:1st derivative}
in Lemma \ref{C-O-derivative} yields that\begin{align*}
 & \mathbb{E}\left[\nabla_{x_{j}}\left(\varphi\left(x_{T}\right)\right)\, J_{0,T}^{k}\right]\\
 & =\mathbb{E}\left[\int_{0}^{T}\!\!\int_{\mathbb{R}_{0}^{m}}\nabla_{x_{j}}\left(\mathfrak{B}_{z}u\left(s,x_{s}\right)\right)\sum_{i=1}^{m}\partial_{z_{i}}\left(g\left(z\right)v_{s}^{ik}\left(\varepsilon,z\right)\right)dz\, ds\right]\\
 & =-\mathbb{E}\left[\int_{0}^{T}\!\!\int_{\mathbb{R}_{0}^{m}}\sum_{i=1}^{m}\partial_{z_{i}}\nabla_{x_{j}}\left(u\left(s,x_{s}+b_{z}\left(\varepsilon,x_{s}\right)\right)\right)v_{s}^{ik}\left(\varepsilon,z\right)d\hat{\mu}\right]\\
 & =-\mathbb{E}\left[\int_{0}^{T}\!\!\int_{\mathbb{R}_{0}^{m}}\sum_{i=1}^{m}\sum_{\beta,\delta=1}^{d}\nabla_{\delta}\nabla_{\beta}u\left(s,x_{s}+b_{z}\left(\varepsilon,x_{s}\right)\right)\left[\left(I_{d}+\nabla b_{z}\left(\varepsilon,x_{s}\right)\right)Z_{s}\right]_{\beta j}\partial_{z_{i}}\left(b_{z}^{\delta}\left(\varepsilon,x_{s}\right)\right)v_{s}^{ik}\left(\varepsilon,z\right)d\hat{\mu}\right]\\
 & \quad-\mathbb{E}\left[\int_{0}^{T}\!\!\int_{\mathbb{R}_{0}^{m}}\sum_{i=1}^{m}\sum_{\beta=1}^{d}\nabla_{\beta}u\left(s,x_{s}+b_{z}\left(\varepsilon,x_{s}\right)\right)\nabla_{x_{j}}\partial_{z_{i}}\left(b_{z}^{\beta}\left(\varepsilon,x_{s}\right)\right)v_{s}^{ik}\left(\varepsilon,z\right)d\hat{\mu}\right]\\
 & =-\mathbb{E}\left[\int_{0}^{T}\!\!\int_{\mathbb{R}_{0}^{m}}\sum_{\beta,\delta=1}^{d}\nabla_{\delta}\nabla_{\beta}u\left(s,x_{s}+b_{z}\left(\varepsilon,x_{s}\right)\right)\left[\left(I_{d}+\nabla b_{z}\left(\varepsilon,x_{s}\right)\right)Z_{s}\right]_{\beta j}\left[\left(I_{d}+\nabla b_{z}\left(\varepsilon,x_{s}\right)\right)Z_{s}\right]_{\delta k}\left|z\right|^{2}d\hat{\mu}\right]\\
 & \quad-\mathbb{E}\left[\int_{0}^{T}\!\!\int_{\mathbb{R}_{0}^{m}}\sum_{i=1}^{m}\sum_{\beta=1}^{d}\nabla_{\beta}u\left(s,x_{s}+b_{z}\left(\varepsilon,x_{s}\right)\right)\nabla_{x_{j}}\partial_{z_{i}}\left(b_{z}^{\beta}\left(\varepsilon,x_{s}\right)\right)v_{s}^{ik}\left(\varepsilon,z\right)d\hat{\mu}\right]\\
 & =-\mathbb{E}\left[\int_{0}^{T}\!\!\int_{\mathbb{R}_{0}^{m}}\nabla_{x_{j}}\nabla_{x_{k}}\left(u\left(s,x_{s}+b_{z}\left(\varepsilon,x_{s}\right)\right)\right)\left|z\right|^{2}d\hat{\mu}\right]\\
 & \quad+\mathbb{E}\left[\int_{0}^{T}\!\!\int_{\mathbb{R}_{0}^{m}}\sum_{\beta=1}^{d}\nabla_{\beta}u\left(s,x_{s}+b_{z}\left(\varepsilon,x_{s}\right)\right)\,\nabla_{x_{k}}\left(\left[\left(I_{d}+\nabla b_{z}\left(\varepsilon,x_{s}\right)\right)Z_{s}\right]_{\beta j}\right)\left|z\right|^{2}d\hat{\mu}\right]\\
 & \quad-\mathbb{E}\left[\int_{0}^{T}\!\!\int_{\mathbb{R}_{0}^{m}}\sum_{i=1}^{m}\sum_{\beta=1}^{d}\nabla_{\beta}u\left(s,x_{s}+b_{z}\left(\varepsilon,x_{s}\right)\right)\nabla_{x_{j}}\partial_{z_{i}}\left(b_{z}^{\beta}\left(\varepsilon,x_{s}\right)\right)v_{s}^{ik}\left(\varepsilon,z\right)d\hat{\mu}\right],\end{align*}
 where Assumption \ref{Levy measure} is used for the second equality. 

Replace $T$ by $\tilde{T}$ and $\varphi$ by $P_{\tilde{T}}\varphi$,
respectively. In a similar manner to Theorem \ref{Theorem 1}, we
have \[
\mathbb{E}\left[\nabla_{x_{j}}\left(P_{\tilde{T}}\varphi\left(x_{\tilde{T}}\right)\right)\, J_{0,\tilde{T}}^{k}\right]=\mathbb{E}\left[\nabla_{x_{j}}\left(\mathbb{E}\left[\varphi\left(x_{T}\right)\big|\mathcal{F}_{\tilde{T}}\right]\right)J_{0,\tilde{T}}^{k}\right]=\mathbb{E}\left[\varphi\left(x_{T}\right)\Gamma_{\tilde{T},T}^{\left(1\right),j}J_{0,\tilde{T}}^{k}\right].\]
Since $\tilde{u}\left(t,x\right)=\mathbb{E}\left[\left(P_{\tilde{T}}\varphi\right)\left(x_{\tilde{T}}\right)\big|x_{0}=x\right]$
for $t\in[0,\tilde{T}]$ and $x\in\mathbb{R}^{d}$, we see that\begin{align*}
 & \mathbb{E}\left[\int_{0}^{\tilde{T}}\!\!\int_{\mathbb{R}_{0}^{m}}\nabla_{x_{j}}\nabla_{x_{k}}\left(\tilde{u}\left(s,x_{s}+b_{z}\left(\varepsilon,x_{s}\right)\right)\right)\left|z\right|^{2}d\hat{\mu}\right]\\
 & =\nabla_{x_{j}}\nabla_{x_{k}}\left(\mathbb{E}\left[\int_{0}^{\tilde{T}}\!\!\int_{\mathbb{R}_{0}^{m}}\tilde{u}\left(s,x_{s}+b_{z}\left(x_{s}\right)\right)\left|z\right|^{2}d\hat{\mu}\right]\right)\\
 & =\nabla_{x_{j}}\nabla_{x_{k}}\left(\mathbb{E}\left[P_{\tilde{T}}\varphi\left(x_{\tilde{T}}\right)\int_{0}^{\tilde{T}}\!\!\int_{\mathbb{R}_{0}^{m}}\left|z\right|^{2}d\mu\right]\right)\\
 & =\nabla_{x_{j}}\nabla_{x_{k}}\left(\mathbb{E}\left[\varphi\left(x_{T}\right)\int_{0}^{\tilde{T}}\!\!\int_{\mathbb{R}_{0}^{m}}\left|z\right|^{2}d\mu\right]\right)\end{align*}
from $\varphi\in C_{K}^{2}\left(\mathbb{R}^{d}\,;\,\mathbb{R}\right)$
and Lemma \ref{key eq}. On the other hand, Assumption \ref{Levy measure}
implies that\begin{align*}
 & \mathbb{E}\left[\int_{0}^{\tilde{T}}\!\!\int_{\mathbb{R}_{0}^{m}}\sum_{\beta=1}^{d}\nabla_{\beta}\tilde{u}\left(s,x_{s}+b_{z}\left(\varepsilon,x_{s}\right)\right)\,\nabla_{x_{k}}\left(\left[\left(I_{d}+\nabla b_{z}\left(\varepsilon,x_{s}\right)\right)Z_{s}\right]_{\beta j}\right)\left|z\right|^{2}d\hat{\mu}\right]\\
 & =-\mathbb{E}\left[\int_{0}^{\tilde{T}}\!\!\int_{\mathbb{R}_{0}^{m}}\sum_{i=1}^{m}\sum_{\beta=1}^{d}\nabla_{\beta}\tilde{u}\left(s,x_{s}+b_{z}\left(\varepsilon,x_{s}\right)\right)\partial_{z_{i}}\left(b_{z}^{\beta}\left(\varepsilon,x_{s}\right)\right)F_{ijk}^{2}\left(\varepsilon,s,z\right)d\hat{\mu}\right]\\
 & =\mathbb{E}\left[\int_{0}^{\tilde{T}}\!\!\int_{\mathbb{R}_{0}^{m}}\mathfrak{B}_{z}\tilde{u}\left(s,x_{s}\right)\,\left\{ \mbox{div}_{z}\left[F^{2}\left(\varepsilon,s,z\right)\right]\right\} _{jk}d\hat{\mu}\right]\\
 & =\mathbb{E}\left[P_{\tilde{T}}\varphi\left(x_{\tilde{T}}\right)\int_{0}^{\tilde{T}}\!\!\int_{\mathbb{R}_{0}^{m}}\left\{ \mbox{div}_{z}\left[F^{2}\left(\varepsilon,s,z\right)\right]\right\} _{jk}d\tilde{\mu}\right]\\
 & =\mathbb{E}\left[\varphi\left(x_{T}\right)\int_{0}^{\tilde{T}}\!\!\int_{\mathbb{R}_{0}^{m}}\left\{ \mbox{div}_{z}\left[F^{2}\left(\varepsilon,s,z\right)\right]\right\} _{jk}d\tilde{\mu}\right]\end{align*}
from \eqref{eq:semi-group} in the proof of Lemma \ref{Continuous 2nd}.
Similarly, we have\begin{align*}
 & -\mathbb{E}\left[\int_{0}^{\tilde{T}}\!\!\int_{\mathbb{R}_{0}^{m}}\sum_{i=1}^{m}\sum_{\beta=1}^{d}\nabla_{\beta}\tilde{u}\left(s,x_{s}+b_{z}\left(\varepsilon,x_{s}\right)\right)\nabla_{x_{j}}\partial_{z_{i}}\left(b_{z}^{\beta}\left(\varepsilon,x_{s}\right)\right)v_{s}^{ik}\left(\varepsilon,z\right)d\hat{\mu}\right]\\
 & =-\mathbb{E}\left[\int_{0}^{\tilde{T}}\!\!\int_{\mathbb{R}_{0}^{m}}\sum_{i=1}^{m}\sum_{\beta=1}^{d}\nabla_{\beta}\tilde{u}\left(s,x_{s}+b_{z}\left(\varepsilon,x_{s}\right)\right)\partial_{z_{i}}\left(b_{z}^{\beta}\left(\varepsilon,x_{s}\right)\right)F_{ijk}^{3}\left(\varepsilon,s,z\right)d\hat{\mu}\right]\\
 & =\mathbb{E}\left[\int_{0}^{\tilde{T}}\!\!\int_{\mathbb{R}_{0}^{m}}\mathfrak{B}_{z}\tilde{u}\left(s,x_{s}\right)\left\{ \mbox{div}_{z}\left[F^{3}\left(\varepsilon,s,z\right)\right]\right\} _{jk}d\hat{\mu}\right]\\
 & =\mathbb{E}\left[P_{\tilde{T}}\varphi\left(x_{\tilde{T}}\right)\int_{0}^{\tilde{T}}\!\!\int_{\mathbb{R}_{0}^{m}}\left\{ \mbox{div}_{z}\left[F^{3}\left(\varepsilon,s,z\right)\right]\right\} _{jk}d\tilde{\mu}\right]\\
 & =\mathbb{E}\left[\varphi\left(x_{T}\right)\int_{0}^{\tilde{T}}\!\!\int_{\mathbb{R}_{0}^{m}}\left\{ \mbox{div}_{z}\left[F^{3}\left(\varepsilon,s,z\right)\right]\right\} _{jk}d\tilde{\mu}\right]\end{align*}
from Assumption \ref{Levy measure} and \eqref{eq:semi-group} in
the proof of Lemma \ref{Continuous 2nd}. The proof of Lemma \ref{Jump 2nd}
is complete. \hfill{}$\square$

\bigskip{}

\begin{cor}\label{Corollary 2}For $1\le j,\ k\le d$ and $\varphi\in C_{K}^{2}\left(\mathbb{R}^{d}\,;\,\mathbb{R}\right)$,
it holds that \begin{equation}
\begin{split} & \nabla_{x_{j}}\nabla_{x_{k}}\left(\mathbb{E}\left[\varphi\left(x_{T}\right)A_{0,\tilde{T}}\right]\right)\\
 & =\mathbb{E}\left[\varphi\left(x_{T}\right)\Gamma_{\tilde{T},T}^{\left(1\right),j}\left(L_{0,\tilde{T}}^{k}-J_{0,\tilde{T}}^{k}\right)\right]\\
 & \quad\ +\mathbb{E}\left[\varphi\left(x_{T}\right)\left\{ \int_{0}^{\tilde{T}}F^{1}\left(\varepsilon,s\right)dW_{s}+\sum_{\sigma=2}^{3}\int_{0}^{\tilde{T}}\!\!\int_{\mathbb{R}_{0}^{m}}\mbox{div}_{z}\left[F^{\sigma}\left(\varepsilon,s,z\right)\right]d\tilde{\mu}\right\} _{jk}\right].\end{split}
\label{eq:2nd derivative}\end{equation}
 \end{cor}

\bigskip{}

\emph{Proof of Theorem \ref{Theorem 2}.} Our goal is to remove $A_{0,\tilde{T}}$
in the left hand side of the equality in Corollary \ref{Corollary 2}.
In order to do that, we shall adopt the same strategy as in Theorem
\ref{Theorem 1}. Define a new probability measure $\tilde{\mathbb{P}}^{\lambda}$
by \[
\frac{d\tilde{\mathbb{P}}^{\lambda}}{d\mathbb{P}}\big|_{\mathcal{F}_{T}}=\exp\left\{ -\lambda A_{0,\tilde{T}}-N_{0,\tilde{T}}^{\lambda}\right\} \]
via the Girsanov transform, and $\tilde{\mathbb{E}}^{\lambda}\left[\ \cdot\ \right]$
is the expectation with respect to the measure $\tilde{\mathbb{P}}^{\lambda}$,
where $N_{\tau,t}^{\lambda}=\int_{\tau}^{t}\!\!\int_{\mathbb{R}_{0}^{m}}\left(e^{-\lambda\left|z\right|^{2}}-1\right)d\hat{\mu}$.
As stated in the proof of Theorem \ref{Theorem 1}, we have \[
\mathbb{E}\left[\varphi\left(x_{T}\right)\right]=\mathbb{E}\left[\frac{\varphi\left(x_{T}\right)A_{0,\tilde{T}}}{A_{0,\tilde{T}}}\right]=\int_{0}^{\infty}\tilde{\mathbb{E}}^{\lambda}\left[\varphi\left(x_{T}\right)A_{0,\tilde{T}}\right]e^{N_{0,\tilde{T}}^{\lambda}}d\lambda,\]
and, under the measure $\tilde{\mathbb{P}}^{\lambda}$, $d\mu$ is
the Poisson random measure with intensity $d\hat{\mu}_{\lambda}:=e^{-\lambda\left|z\right|^{2}}d\hat{\mu}$.
Moreover, $d\tilde{\mu}_{\lambda}:=d\mu-d\hat{\mu}_{\lambda}$ is
a martingale measure (cf. \cite{Kunita}). 

Define $\Gamma_{\tau,t}^{(1),(\lambda)}=\left(L_{\tau,t}-J_{\tau,t}^{(\lambda)}\right)/A_{\tau,t}+K_{\tau,t}/\left(A_{\tau,t}\right)^{2}$.
Applying Corollary \ref{Corollary 2} with respect to $\tilde{\mathbb{E}}^{\lambda}\left[\ \cdot\ \right]$,
we have\begin{align*}
 & \nabla_{x_{j}}\nabla_{x_{k}}\left(\mathbb{E}\left[\varphi\left(x_{T}\right)\right]\right)\\
 & =\int_{0}^{\infty}e^{N_{0,\tilde{T}}^{\lambda}}\nabla_{x_{j}}\nabla_{x_{k}}\left(\tilde{\mathbb{E}}^{\lambda}\left[\varphi\left(x_{T}\right)A_{0,\tilde{T}}\right]\right)d\lambda\\
 & =\int_{0}^{\infty}e^{N_{0,\tilde{T}}^{\lambda}}\Bigg\{\tilde{\mathbb{E}}^{\lambda}\left[\varphi\left(x_{T}\right)\Gamma_{\tilde{T},T}^{\left(1\right),(\lambda),j}\left(L_{0,\tilde{T}}^{k}-J_{0,\tilde{T}}^{\left(\lambda\right),k}\right)\right]\\
 & \qquad\qquad\negthickspace+\tilde{\mathbb{E}}^{\lambda}\left[\varphi\left(x_{T}\right)\left\{ \int_{0}^{\tilde{T}}\left\{ F^{1}\left(\varepsilon,s\right)\right\} _{jk}dW_{s}+\sum_{\sigma=2}^{3}\int_{0}^{\tilde{T}}\!\!\int_{\mathbb{R}_{0}^{m}}\left\{ \mbox{div}_{z}\left[F^{\sigma}\left(\varepsilon,s,z\right)\right]\right\} _{jk}d\tilde{\mu}_{\lambda}\right\} \right]\Bigg\}d\lambda\\
 & =:I_{1}+I_{2},\end{align*}
where $J_{\tau,t}^{\left(\lambda\right)}=\int_{\tau}^{t}\!\!\int_{\mathbb{R}_{0}^{m}}\left\{ \mbox{div}_{z}\left[e^{-\lambda\left|z\right|^{2}}g\left(z\right)v_{s}\left(\varepsilon,z\right)\right]\right\} _{k}\Big/\left(e^{-\lambda\left|z\right|^{2}}g\left(z\right)\right)\, d\tilde{\mu}_{\lambda}$
for $0\le\tau\le t\le T$. Since $d\tilde{\mu}_{\lambda}=d\tilde{\mu}+\left(1-e^{-\lambda\left|z\right|^{2}}\right)d\hat{\mu}$
and $J_{\tau,t}^{\left(\lambda\right)}=J_{\tau,t}-\lambda K_{\tau,t}$
as seen in the proof of Theorem \ref{Theorem 1}, we have \[
\Gamma_{\tau,t}^{(1),(\lambda)}=\frac{\, L_{\tau,t}-J_{\tau,t}\,}{A_{\tau,t}}+\frac{K_{\tau,t}}{\left(A_{\tau,t}\right)^{2}}+\lambda\frac{K_{\tau,t}}{\, A_{\tau,t}\,}.\]
Thus, it holds that \[
\begin{split}I_{1} & =\mathbb{E}\left[\varphi\left(x_{T}\right)\int_{0}^{\infty}e^{-\lambda A_{0,\tilde{T}}}\left(\Gamma_{\tilde{T},T}^{(1),j}+\lambda\frac{K_{\tilde{T},T}^{j}}{\, A_{\tilde{T},T}\,}\right)\left(L_{0,\tilde{T}}^{k}-J_{0,\tilde{T}}^{k}+\lambda K_{0,\tilde{T}}^{k}\right)d\lambda\right]\\
 & =\mathbb{E}\left[\varphi\left(x_{T}\right)\int_{0}^{\infty}e^{-\lambda A_{0,\tilde{T}}}\left\{ \Gamma_{\tilde{T},T}^{(1),j}\left(L_{0,\tilde{T}}^{k}-J_{0,\tilde{T}}^{k}+\lambda K_{0,\tilde{T}}^{k}\right)+\lambda\frac{K_{\tilde{T},T}^{j}}{\, A_{\tilde{T},T}\,}\left(L_{0,\tilde{T}}^{k}-J_{0,\tilde{T}}^{k}\right)+\lambda^{2}\frac{K_{\tilde{T},T}^{j}}{\, A_{\tilde{T},T}\,}K_{0,\tilde{T}}^{k}\right\} \, d\lambda\right]\\
 & =\mathbb{E}\left[\varphi\left(x_{T}\right)\left\{ \Gamma_{\tilde{T},T}^{(1),j}\,\Gamma_{0,\tilde{T}}^{(1),k}+\frac{K_{\tilde{T},T}^{j}}{\, A_{\tilde{T},T}\,}\,\frac{\, L_{0,\tilde{T}}^{k}-J_{0,\tilde{T}}^{k}\,}{\left(A_{0,\tilde{T}}\right)^{2}}+2\frac{K_{\tilde{T},T}^{j}}{\, A_{\tilde{T},T}\,}\,\frac{K_{0,\tilde{T}}^{k}}{\,\left(A_{0,\tilde{T}}\right)^{3}\,}\right\} \right]\end{split}
\]
from the Fubini theorem. Similarly, we have\begin{align*}
I_{2} & =\mathbb{E}\left[\varphi\left(x_{T}\right)\int_{0}^{\infty}e^{-\lambda A_{0,\tilde{T}}}\left\{ \int_{0}^{\tilde{T}}F^{1}\left(\varepsilon,s\right)\, dW_{s}+\int_{0}^{\tilde{T}}\!\!\int_{\mathbb{R}_{0}^{m}}\sum_{\sigma=2}^{3}\mbox{div}_{z}\left[F^{\sigma}\left(\varepsilon,s,z\right)\right]\, d\tilde{\mu_{\lambda}}\right\} _{jk}d\lambda\right]\\
 & =\mathbb{E}\left[\varphi\left(x_{T}\right)\,\frac{1}{\, A_{0,\tilde{T}}\,}\left\{ \int_{0}^{\tilde{T}}F^{1}\left(\varepsilon,s\right)\, dW_{s}+\int_{0}^{\tilde{T}}\!\!\int_{\mathbb{R}_{0}^{m}}\sum_{\sigma=2}^{3}\mbox{div}_{z}\left[F^{\sigma}\left(\varepsilon,s,z\right)\right]d\mu\right\} _{jk}\right]\\
 & \qquad-\mathbb{E}\left[\varphi\left(x_{T}\right)\left\{ \sum_{\sigma=2}^{3}\int_{0}^{\tilde{T}}\!\!\int_{\mathbb{R}_{0}^{m}}\frac{\mbox{div}_{z}\left[F^{\sigma}\left(\varepsilon,s,z\right)\right]}{A_{0,\tilde{T}}+\left|z\right|^{2}}d\hat{\mu}\right\} _{jk}\right].\end{align*}
The proof of Theorem \ref{Theorem 2} is complete. \hfill{}$\square$

\subsection{Proof of Corollary \ref{sensitivity formulae}}

For $\varphi\in C_{K}^{2}\left(\mathbb{R}^{d}\,;\,\mathbb{R}\right)$,
all sensitivity formulae are the direct consequences of Theorem \ref{Theorem 1},
\ref{vega} and \ref{Theorem 2}. The strategy to remove the regularity
conditions, and to extend to the class $\mathfrak{F}$, is almost
parallel to the one studied in \cite{Kawai-Takeuchi}. 

First, we shall extend from $C_{K}^{2}\left(\mathbb{R}^{d}\,;\,\mathbb{R}\right)$
to $C_{K}\left(\mathbb{R}^{d}\,;\,\mathbb{R}\right)$. Since $\varphi\in C_{K}\left(\mathbb{R}^{d}\,;\,\mathbb{R}\right)$
can be approximated uniformly and boundedly by a sequence $\left\{ \varphi_{n}\,;\, n\in\mathbb{N}\right\} $,
we see that, for each compact set $H\subset\mathbb{R}^{d}$, \begin{gather*}
\left|\mathbb{E}\left[\varphi\left(x_{T}\right)\right]-\mathbb{E}\left[\varphi_{n}\left(x_{T}\right)\right]\right|\le\left\Vert \varphi-\varphi_{n}\right\Vert _{\infty},\\
\sup_{x\in H}\left|\nabla_{x}\left(\mathbb{E}\left[\varphi_{n}\left(x_{T}\right)\right]\right)-\mathbb{E}\left[\varphi\left(x_{T}\right)\,\Gamma_{T}^{\left(1\right)}\right]\right|\le\sup_{x\in H}\mathbb{E}\left[\left|\Gamma_{T}^{\left(1\right)}\right|^{2}\right]^{1/2}\left\Vert \varphi_{n}-\varphi\right\Vert _{\infty},\end{gather*}
which tends to $0$ as $n\to+\infty$. Thus, the sensitivity formula
$\nabla_{x}\left(\mathbb{E}\left[\varphi\left(x_{T}\right)\right]\right)=\mathbb{E}\left[\varphi\left(x_{T}\right)\Gamma_{T}^{\left(1\right)}\right]$
holds for $\varphi\in C_{K}\left(\mathbb{R}^{d}\,;\,\mathbb{R}\right)$. 

Second, we shall extend to the class $C_{b}\left(\mathbb{R}^{d}\,;\,\mathbb{R}\right)$
of bounded continuous functions. Let $\sigma\in\left(0,1\right)$
be fixed, and write $N\left(y;\delta\right)=\left\{ \tilde{y}\in\mathbb{R}^{d}\,;\,\left|\tilde{y}-y\right|<\delta\right\} $
for $y\in\mathbb{R}^{d}$ and $\delta>0$. For $\varphi\in C_{b}\left(\mathbb{R}^{d}\,;\,\mathbb{R}\right)$,
we can find the sequence $\left\{ \varphi_{n}\,;\, n\in\mathbb{N}\right\} $
of continuous functions defined by \[
\varphi_{n}\left(x\right)=\begin{cases}
\varphi\left(x\right), & \mbox{if }x\in\overline{N\left(\bm{0};n-\sigma\right)},\\
0, & \mbox{if }x\in N\left(\bm{0};n+\sigma\right)^{c},\end{cases}\]
and $\varphi_{n}\left(x\right)\in\left[0,\varphi\left(x\right)\right]$
for each $x\in\left(\overline{N\left(\bm{0};n-\sigma\right)}\right)^{c}\cap N\left(\bm{0};n+\sigma\right)$,
where $\left[0,-1\right]$ should be understood as $\left[-1,0\right]$.
Clearly, $\varphi_{n}\in C_{K}\left(\mathbb{R}^{d}\,;\,\mathbb{R}\right)$,
and $\sup_{n\in\mathbb{N}}\left\Vert \varphi_{n}\right\Vert _{\infty}=\left\Vert \varphi\right\Vert _{\infty}$.
The dominated convergence theorem leads to \[
\left|\mathbb{E}\left[\varphi\left(x_{T}\right)\right]-\mathbb{E}\left[\varphi_{n}\left(x_{T}\right)\right]\right|\to0\]
as $n\to+\infty$. On the other hand, since \[
\mathbb{E}\left[\left|\varphi_{n}\left(x_{T}\right)-\varphi\left(x_{T}\right)\right|^{2}\right]\le\mathbb{E}\left[\left|\varphi\left(x_{T}\right)\right|^{2}I_{\left(\left|x_{T}\right|>n-\sigma\right)}\right]\le\frac{\left\Vert \varphi\right\Vert _{\infty}^{2}}{\left(n-\sigma\right)^{2}}\mathbb{E}\left[\left|x_{T}\right|^{2}\right]\to0\]
as $n\to+\infty$, we have \[
\sup_{x\in H}\left|\nabla_{x}\left(\mathbb{E}\left[\varphi_{n}\left(x_{T}\right)\right]\right)-\mathbb{E}\left[\varphi\left(x_{T}\right)\,\Gamma_{T}^{\left(1\right)}\right]\right|\le\sup_{x\in H}\mathbb{E}\left[\left|\Gamma_{T}^{\left(1\right)}\right|^{2}\right]^{1/2}\sup_{x\in H}\mathbb{E}\left[\left|\varphi_{n}\left(x_{T}\right)-\varphi\left(x_{T}\right)\right|^{2}\right]^{1/2},\]
which tends to $0$ as $n\to+\infty$. Hence, we can obtain the sensitivity
formula $\nabla_{x}\left(\mathbb{E}\left[\varphi\left(x_{T}\right)\right]\right)=\mathbb{E}\left[\varphi\left(x_{T}\right)\Gamma_{T}^{\left(1\right)}\right]$
for $\varphi\in C_{b}\left(\mathbb{R}^{d}\,;\,\mathbb{R}\right)$.

Thirdly, we shall extend to the class of finite linear combinations
of indicator functions, which leads us to extend to the class $\mathfrak{F}$
immediately, via the standard truncation argument. It is sufficient
to consider the case $\varphi=I_{U}$ for a subset $U$ in $\mathbb{R}^{d}$.
Then, we can find a sequence $\left\{ \varphi_{n}\,;\, n\in\mathbb{N}\right\} $
of continuous functions such that \[
\varphi_{n}\left(x\right)=\begin{cases}
\varphi\left(x\right), & \mbox{if }x\in U_{-,}\\
0, & \mbox{if }x\in U_{+}^{c},\end{cases}\]
and $\varphi_{n}\left(x\right)\in\left[0,\varphi\left(x\right)\right]$
for $x\in U_{-}^{c}\cap U_{+}$, where \begin{gather*}
U_{+}=\left\{ y\in\mathbb{R}^{d}\,;\,\left|y-\tilde{y}\right|<\frac{1}{n}\ \left(\tilde{y}\in\partial\overline{U}\right)\right\} \cup U,\ U_{-}=\left\{ y\in\mathbb{R}^{d}\,;\,\left|y-\tilde{y}\right|<\frac{1}{n}\ \left(\tilde{y}\in\partial\overline{U}\right)\right\} \cap U.\end{gather*}
Clearly, $\varphi_{n}\in C_{b}\left(\mathbb{R}^{d}\,;\,\mathbb{R}\right)$,
and $\sup_{n\in\mathbb{N}}\left\Vert \varphi_{n}\right\Vert _{\infty}\le1$.
The dominated convergence theorem implies that \[
\left|\mathbb{E}\left[\varphi\left(x_{T}\right)\right]-\mathbb{E}\left[\varphi_{n}\left(x_{T}\right)\right]\right|\to0\]
as $n\to+\infty$. On the other hand, since there exists a smooth
density $p_{T}\left(\varepsilon,x,y\right)$ for the random variable
$x_{T}$ with respect to the Lebesgue measure on $\mathbb{R}^{d}$
as stated in Proposition \ref{known fact}, we have\begin{align*}
\sup_{x\in H}\mathbb{E}\left[\left|\varphi_{n}\left(x_{T}\right)-\varphi\left(x_{T}\right)\right|^{2}\right] & =\sup_{x\in H}\mathbb{E}\left[\left|\varphi_{n}\left(x_{T}\right)-\varphi\left(x_{T}\right)\right|^{2};x_{T}\in U_{-}^{c}\cap U_{+}\right]\\
 & \le4\,\left|U_{-}^{c}\cap U_{+}\right|\,\sup_{x\in H}\sup_{y\in\overline{U_{-}^{c}\cap U_{+}}}p_{T}\left(\varepsilon,x,y\right),\end{align*}
which tends to $0$ as $n\to+\infty$, because of $\left|\overline{U_{-}^{c}\cap U_{+}}\right|\to0$.
Hence, we have \[
\sup_{x\in H}\left|\nabla_{x}\left(\mathbb{E}\left[\varphi_{n}\left(x_{T}\right)\right]\right)-\mathbb{E}\left[\varphi\left(x_{T}\right)\,\Gamma_{T}^{\left(1\right)}\right]\right|\le\sup_{x\in H}\mathbb{E}\left[\left|\Gamma_{T}^{\left(1\right)}\right|^{2}\right]^{1/2}\sup_{x\in H}\mathbb{E}\left[\left|\varphi_{n}\left(x_{T}\right)-\varphi\left(x_{T}\right)\right|^{2}\right]^{1/2},\]
which tends to $0$ as $n\to+\infty$. Therefore, we can conclude
that the sensitivity formula \[
\nabla_{x}\left(\mathbb{E}\left[\varphi\left(x_{T}\right)\right]\right)=\mathbb{E}\left[\varphi\left(x_{T}\right)\Gamma_{T}^{\left(1\right)}\right]\]
 holds for $\varphi\in\mathfrak{F}$.

The regularity condition on the function $\varphi$ in Theorem \ref{vega}
and \ref{Theorem 2} can be relaxed to the class $\mathfrak{F}$ in
a similar manner. \hfill{}$\square$

\section{Examples}

\begin{exa}[L\'evy processes]\rm{\label{Levy}Let $m=d=1$, and
$\left(x,\gamma,\sigma_{1},\sigma_{2}\right)\in\mathbb{R}^{4}$. Consider
the $\mathbb{R}$-valued process $\left\{ x_{t}\,;\, t\in\left[0,T\right]\right\} $
given by \[
x_{t}=x+\gamma t+\sigma_{1}W_{t}+\sigma_{2}\int_{0}^{t}\!\!\int_{\mathbb{R}_{0}}z\, d\overline{\mu}.\]

Consider the case of $\sigma_{1}\neq0$ and $\sigma_{2}\neq0$. Since
\begin{gather*}
A_{\tau,t}=\left(t-\tau\right)+\int_{\tau}^{t}\!\!\int_{\mathbb{R}_{0}}\left|z\right|^{2}d\mu,\quad L_{\tau,t}=\frac{W_{t}-W_{\tau}}{\sigma_{1}},\quad J_{\tau,t}=\int_{\tau}^{t}\!\!\int_{\mathbb{R}_{0}}\frac{\left\{ g\left(z\right)\,\left|z\right|^{2}\right\} ^{\prime}}{\sigma_{2}g\left(z\right)}\, d\tilde{\mu},\ \\
K_{\tau,t}=\int_{\tau}^{t}\!\!\int_{\mathbb{R}_{0}}\frac{z}{\sigma_{2}}\, d\mu,\quad L_{T}^{\gamma}=\frac{W_{T}}{\sigma_{1}},\quad R_{T}^{\sigma_{1}}-Q_{T}^{\sigma_{1}}=\frac{W_{T}^{2}-T}{\sigma_{1}T},\ \\
L_{T}^{\sigma_{2}}=-\frac{W_{T}}{\sigma_{1}}\int_{\left|z\right|\le1}z\, d\nu,\quad J_{T}^{\sigma_{2}}=\int_{0}^{T}\!\!\int_{\mathbb{R}_{0}}\frac{\left\{ g\left(z\right)\,\left|z\right|^{2}\right\} ^{\prime}}{\sigma_{2}\, g\left(z\right)}d\tilde{\mu},\end{gather*}
we have \begin{align*}
\Gamma_{T}^{\left(1\right)} & =\frac{L_{0,T}-J_{0,T}}{A_{0,T}}+\frac{K_{0,T}}{\left(A_{0,T}\right)^{2}},\quad\Gamma_{T}^{\left(2,\gamma\right)}=L_{T}^{\gamma},\quad\Gamma_{T}^{\left(2,\sigma_{1}\right)}=R_{T}^{\sigma_{1}}-Q_{T}^{\sigma_{1}},\quad\Gamma_{T}^{\left(2,\sigma_{2}\right)}=L_{T}^{\sigma_{2}}+J_{T}^{\sigma_{2}},\\
\Gamma_{T}^{\left(3\right)} & =\left\{ \Gamma_{\tilde{T},T}^{(1)}+\frac{K_{\tilde{T},T}}{A_{\tilde{T},T}\, A_{0,\tilde{T}}}\right\} \,\Gamma_{0,\tilde{T}}^{(1)}+\frac{K_{\tilde{T},T}\, K_{0,\tilde{T}}}{A_{\tilde{T},T}\,\left(A_{0,\tilde{T}}\right)^{3}}.\end{align*}

As stated in Remark \ref{uniformly elliptic diffusion} and \ref{uniformly elliptic jump},
the case of either $\sigma_{1}\neq0$ or $\sigma_{2}\neq0$ is also
in our position. In the case of $\sigma_{1}\neq0$, since \[
L_{\tau,t}=\frac{W_{t}-W_{\tau}}{\sigma_{1}},\quad L_{T}^{\gamma}=\frac{W_{T}}{\sigma_{1}},\quad R_{T}^{\sigma_{1}}-Q_{T}^{\sigma_{1}}=\frac{W_{T}^{2}-T}{\sigma_{1}T},\quad L_{T}^{\sigma_{2}}=-\frac{W_{T}}{\sigma_{1}}\int_{\left|z\right|\le1}z\, d\nu,\]
we have\begin{gather*}
\Gamma_{T}^{\left(1\right)}=\frac{L_{0,T}}{T},\quad\Gamma_{T}^{\left(2,\gamma\right)}=L_{T}^{\gamma},\quad\Gamma_{T}^{\left(2,\sigma_{1}\right)}=R_{T}^{\sigma_{1}}-Q_{T}^{\sigma_{1}},\quad\Gamma_{T}^{\left(2,\sigma_{2}\right)}=L_{T}^{\sigma_{2}},\quad\Gamma_{T}^{\left(3\right)}=\frac{L_{0,\tilde{T}}\, L_{\tilde{T},T}}{\tilde{T}^{2}}.\end{gather*}
In the case of $\sigma_{2}\neq0$, since\begin{gather*}
J_{\tau,t}=\int_{\tau}^{t}\!\!\int_{\mathbb{R}_{0}}\frac{\left\{ g\left(z\right)\,\left|z\right|^{2}\right\} ^{\prime}}{\sigma_{2}g\left(z\right)}\, d\tilde{\mu},\quad K_{\tau,t}=\int_{\tau}^{t}\!\!\int_{\mathbb{R}_{0}}\frac{z}{\sigma_{2}}\, d\mu,\quad J_{T}^{\sigma_{2}}=\int_{0}^{T}\!\!\int_{\mathbb{R}_{0}}\frac{\left\{ g\left(z\right)\,\left|z\right|^{2}\right\} ^{\prime}}{\sigma_{2}\, g\left(z\right)}d\tilde{\mu},\end{gather*}
 we have\begin{gather*}
\Gamma_{T}^{\left(1\right)}=-\frac{V_{0,T}}{A_{0,T}}+\frac{K_{0,T}}{\left(A_{0,T}\right)^{2}},\quad\Gamma_{T}^{\left(2,\sigma_{2}\right)}=J_{T}^{\sigma_{2}},\quad\Gamma_{T}^{\left(3\right)}=\left\{ \Gamma_{\tilde{T},T}^{(1)}+\frac{K_{\tilde{T},T}}{A_{\tilde{T},T}\, A_{0,\tilde{T}}}\right\} \,\Gamma_{0,\tilde{T}}^{(1)}+\frac{K_{\tilde{T},T}\, K_{0,\tilde{T}}}{A_{\tilde{T},T}\,\left(A_{0,\tilde{T}}\right)^{3}},\end{gather*}
where $A_{\tau,t}=\int_{\tau}^{t}\!\!\int_{\mathbb{R}_{0}}\left|z\right|^{2}d\mu$.
\hfill{}$\square$

}\end{exa}\bigskip{}

\begin{exa}[geometric L\'evy processes]\rm{Let $m=d=1$, $\left(\gamma,\sigma_{1},\sigma_{2}\right)\in\mathbb{R}^{3},$
and $\left\{ X_{t}\,;\, t\in\left[0,T\right]\right\} $ the $\mathbb{R}$-valued
L\'evy process represented as follows: \[
X_{t}=\gamma t+\sigma_{1}W_{t}+\sigma_{2}\int_{0}^{t}\!\!\int_{\mathbb{R}_{0}}z\, d\overline{\mu}.\]
Let $x>0$, and $\left\{ x_{t}\,;\, t\in\left[0,T\right]\right\} $
the $\mathbb{R}$-valued process defined by $x_{t}=xe^{X_{t}}$, which
is called \textit{the geometric L\'evy process}. Let $\varphi\in\mathfrak{F}$
be bounded. Since\begin{align*}
\nabla_{x}\left(\mathbb{E}\left[\varphi\left(xe^{X_{t}}\right)\right]\right) & =\mathbb{E}\left[\varphi^{\prime}\left(xe^{X_{t}}\right)\, e^{X_{t}}\right]=\frac{1}{x}\nabla_{X}\left(\mathbb{E}\left[\varphi\left(e^{X+X_{t}}\right)\right]\right)\big|_{X=\log x},\\
\partial_{\gamma}\left(\mathbb{E}\left[\varphi\left(x_{T}\right)\right]\right) & =\mathbb{E}\left[\varphi^{\prime}\left(x_{T}\right)\, x_{T}\, T\right]=\partial_{\gamma}\left(\mathbb{E}\left[\left(\varphi\circ\psi\right)\left(X+X_{T}\right)\right]\right)\big|_{X=\log x},\\
\partial_{\sigma_{1}}\left(\mathbb{E}\left[\varphi\left(x_{T}\right)\right]\right) & =\mathbb{E}\left[\varphi^{\prime}\left(x_{T}\right)\, x_{T}\, W_{T}\right]=\partial_{\sigma_{1}}\left(\mathbb{E}\left[\left(\varphi\circ\psi\right)\left(X+X_{T}\right)\right]\right)\big|_{X=\log x},\\
\partial_{\sigma_{2}}\left(\mathbb{E}\left[\varphi(x_{T})\right]\right) & =\mathbb{E}\left[\varphi^{\prime}(x_{T})\, x_{T}\,\int_{0}^{T}\!\!\int_{\mathbb{R}_{0}}z\, d\overline{\mu}\right]=\partial_{\sigma_{2}}\left(\mathbb{E}\left[\left(\varphi\circ\psi\right)\left(X+X_{T}\right)\right]\right)\big|_{X=\log x},\\
\nabla_{x}^{2}\left(\mathbb{E}\left[\varphi\left(xe^{X_{t}}\right)\right]\right) & =\mathbb{E}\left[\varphi^{\prime\prime}\left(xe^{X_{t}}\right)\, e^{2X_{t}}\right]=\frac{1}{x^{2}}\left\{ \nabla_{X}^{2}\left(\mathbb{E}\left[\varphi\left(e^{X+X_{t}}\right)\right]\right)-\nabla_{X}\left(\mathbb{E}\left[\varphi\left(e^{X+X_{t}}\right)\right]\right)\right\} \big|_{X=\log x},\end{align*}
 we can calculate the corresponding weights $\tilde{\Gamma}_{T}^{\left(1\right)},\,\tilde{\Gamma}_{T}^{\left(2,\gamma\right)},\,\tilde{\Gamma}_{T}^{\left(2,\sigma_{1}\right)},\,\tilde{\Gamma}_{T}^{\left(2,\sigma_{2}\right)}$
and $\tilde{\Gamma}_{T}^{\left(3\right)}$ by using the result in
Example \ref{Levy}. \hfill{}$\square$

}\end{exa}\bigskip{}

}}
\end{document}